\documentclass[a4paper,12pt]{amsart}
\usepackage{amssymb,a4wide,graphicx,stmaryrd,fullpage,setspace,microtype,textcomp,tikz,enumitem,accents,arydshln,tikz-cd}
\definecolor{dark-red}{rgb}{0.5,0.15,0.15}
\usepackage{hyperref}
\hypersetup{
	colorlinks   = true, 
	urlcolor     = dark-red, 
	linkcolor    = black, 
	citecolor   = dark-red 
}
\usepackage{float}
\usepackage[T1]{fontenc}
\usepackage{lmodern}
\usepackage[numbers,sort&compress]{natbib}

\title{Homotopy theory of Moore flows (III)}

\author[P. Gaucher]{Philippe Gaucher}

\address{Universit\'e Paris Cit\'e, CNRS, IRIF, F-75013, Paris, France}

\urladdr{https://www.irif.fr/{\~{}}gaucher}

\makeatletter
\@namedef{subjclassname@2020}{%
	\textup{2020} Mathematics Subject Classification}
\makeatother
\subjclass[2020]{18C35,18D20,55U35,68Q85}

\keywords{directed path, reparametrization, enriched semicategory, semimonoidal structure, combinatorial model category, Quillen equivalence, locally presentable category, topologically enriched category}

\newcommand{\C}{\mathcal{C}}
\newcommand{\D}{\mathcal{D}}
\newcommand{\K}{\mathcal{K}}
\newcommand{\de}{\partial}
\newcommand{\p}{\times}
\newcommand{\PA}{\mathbb{P}}

\makeatletter
\let\P\@undefined
\let\leq\@undefined
\let\top\@undefined
\makeatother

\newtheorem*{thmN}{Theorem}
\newtheorem{thm}{Theorem}
\newtheorem{prop}{Proposition}
\newtheorem{lem}{Lemma}
\newtheorem{cor}{Corollary}
\newcommand{\bp}{\begin{prop}}
\newcommand{\ep}{\end{prop}}
\newcommand{\bth}{\begin{thm}}
\renewcommand{\eth}{\end{thm}}
\newcommand{\bpf}{\begin{proof}}
\newcommand{\epf}{\end{proof}}

\theoremstyle{definition}
\newtheorem{defn}{Definition}
\newtheorem{rem}{Remark}

\newcommand{\bd}{\begin{defn}}
\newcommand{\ed}{\end{defn}}
\newtheorem{nota}{Notation}

\newtheorem{conj}{Conjecture}

\newcommand{\topspace}{{\mathbf{Top}}}
\newcommand{\iso}{\cong}
\newcommand{\vI}{\overrightarrow{I}}

\newcommand{\moore}{{\mathbb{M}}}
\newcommand{\lmoore}{\mathbb{M}_!}

\newcommand{\topdgrG}{[\mathcal{G}^{op},\topspace]}
\newcommand{\topdgrM}{[\mathcal{M}^{op},\topspace]}
\newcommand{\topdgrP}{[\mathcal{P}^{op},\topspace]}
\newcommand{\brm}[1]{\mathrm{\mathbf{#1}}}

\newcommand{\dtop}{{\brm{Flow}}}
\newcommand{\dtopM}{{\mathcal{M}\brm{Flow}}}
\newcommand{\dtopG}{{\mathcal{G}\brm{Flow}}}
\newcommand{\dtopP}{{\mathcal{P}\brm{Flow}}}

\newcommand{\ttop}{{\brm{TOP}}}
\newcommand{\mtop}{{\brm{MTop}}}
\newcommand{\glob}{{\mathrm{Glob}}}
\newcommand{\globM}{{\mathrm{Glob}}^{\mathcal{M}}}
\newcommand{\globP}{{\mathrm{Glob}}^{\mathcal{P}}}
\newcommand{\globG}{{\mathrm{Glob}}^{\mathcal{G}}}
\DeclareMathOperator{\id}{Id}
\DeclareMathOperator{\Obj}{Obj}

\DeclareMathOperator{\pr}{pr}
\newcommand{\liminj}{\varinjlim}
\newcommand{\dcat}{{\mathrm{cat}}}
\newcommand{\rest}{\!\upharpoonright\!}
\DeclareMathOperator{\carrier}{Carrier}
\DeclareMathOperator{\lan}{Lan}
\newcommand{\ptop}[1]{{\brm{{#1}dTop}}}
\newcommand{\ptopsat}[1]{{\brm{{#1}dTop^{sat}}}}

\newcommand{\cocartesian}{\arrow[lu, phantom, "\ulcorner"{font=\Large}, pos=0]}

\makeatletter
\newcommand*{\@opargbegintheorem}[3]{\trivlist
	\item[\hskip \labelsep{\bfseries #1\ #2}] \textbf{(#3)}\ \itshape}
\makeatother

\newcommand{\ddownarrow}{{\downarrow}}
\DeclareMathOperator{\cocyl}{{Path}}
\newcommand{\ot}{\otimes}
\newcommand{\ttt}{two-out-of-three property}
\DeclareMathOperator{\natgl}{\underline{nat}^{gl}}
\DeclareMathOperator{\natcub}{\underline{nat}^\square}

\allowdisplaybreaks

\begin{document}

\begin{abstract} 
	The previous paper of this series shows that the q-model categories of $\mathcal{G}$-multipointed $d$-spaces and of $\mathcal{G}$-flows are Quillen equivalent. In this paper, the same result is established by replacing the reparametrization category $\mathcal{G}$ by the reparametrization category $\mathcal{M}$. Unlike the case of $\mathcal{G}$, the execution paths of a cellular $\mathcal{M}$-multipointed $d$-space can have stop intervals. The technical tool to overcome this obstacle is the notion of globular naturalization. It is the globular analogue of Raussen's naturalization of a directed path in the geometric realization of a precubical set. The notion of globular naturalization working both for $\mathcal{G}$ and $\mathcal{M}$, the proof of the Quillen equivalence we obtain is valid for the two reparametrization categories. Together with the results of the first paper of this series, we then deduce that $\mathcal{G}$-multipointed $d$-spaces and $\mathcal{M}$-multipointed $d$-spaces have Quillen equivalent q-model structures. Finally, we prove that the saturation hypothesis can be added without any modification in the main theorems of the paper.
\end{abstract}

\maketitle

\tableofcontents
\hypersetup{linkcolor = dark-red}

\section{Introduction}

\subsection*{Presentation of the paper}

This work is a sequel of \cite{Moore1,Moore2} establishing a zigzag of Quillen equivalences between the q-model structures of multipointed $d$-spaces \cite{mdtop} and of flows \cite{model3} thanks to the notion of \textit{Moore flow}. This paper was not initially planned to be presented as a third part of this series. The reason is that, unexpectedly, all proofs of this paper work for the second paper \cite{Moore2} of this series as well thanks to the discovery of a globular analogue of Raussen's notion of naturalization of a directed path (see below). 

Multipointed $d$-spaces and flows are two \textit{multipointed} geometric models of concurrency. This research belongs to a branch of mathematics sometimes called \textit{directed algebraic topology} (DAT) or \textit{directed homotopy theory} \cite{DAT_book}. Note that the latter terminology is also sometimes used for other mathematical research like rewriting system or directed HoTT. DAT studies the homotopical properties of geometric models of concurrency from various points of view. The general idea is that two directed paths which are homotopy equivalent in an appropriate directed sense represent two non-distinguishable possible execution paths of the corresponding concurrent system. The typical example is the one of a full $n$-cube $[0,1]^n$: each continuous path from the initial state $(0,0,\dots,0)$ to the final state $(1,1,\dots,1)$ of the full $n$-cube which is nondecreasing with respect to each axis of coordinates represents the concurrent execution of $n$ actions. Each axis of coordinates represents one action: $0$ means that it is not started and $1$ that it is finished. Nondecreasingness models time irreversibility. In the case of $[0,1]^n$, all directed paths are homotopy equivalent in a directed sense. More general concurrent systems can be modeled by pasting together cubes of various dimensions. The combinatorial notion of precubical set is adapted for such a purpose \cite{DAT_book,zbMATH07226006}. 

The main problem posed by DAT is that the directed segment is \textit{not always} contractible in a directed sense otherwise the causal information could be lost by the associated weak equivalences. For example, contracting the directed segment going from $A$ to $B$ in the branching $C\leftarrow A \rightarrow B$ removes the nondeterministic branching and therefore changes the causal structure. However, the full $n$-cube $[0,1]^n$ is the same object in DAT as the full $n$-cube $[0,2]^n$, which means that the directed segment can be dilated. Moreover, it is possible in $C\leftarrow A \rightarrow B \rightarrow D$ to contract the directed segment going from $A$ to $B$ without changing the causal structure. The non-conventional behavior of the directed segment is the reason why many model category structures introduced in DAT fail to preserve the causal structure. This does not necessarily mean that they are not interesting, just that they probably need to be modified.

Any  geometric model of concurrency encoding directed paths one way or another, including the \textit{non-multipointed} or \textit{continuous} ones of Grandis' $d$-spaces \cite{mg} or Krishnan's streams \cite{MR2545830}, gives rise to a family of spaces $P^1_{\alpha,\beta}$ of \textit{nonconstant directed paths} (the $1$ meaning here of length $1$; they are also called \textit{execution paths} in the sequel) from $\alpha$ to $\beta$ closed under composition and nondecreasing reparametrization. The points $\alpha$ and $\beta$ belong to some set of states chosen in the underlying topological space for continuous models or they run over the set of states for a multipointed model. Choosing this set of states in a continuous model such that, together with the family of spaces $P^1_{\alpha,\beta}$, the information contained in the causal structure is preserved, is related to the research about component categories \cite{DirectedComponents0,DirectedComponents,Raussen2019,Ziemiaski2018}. The latter aims at reducing the size of the fundamental category, and in particular the size of the state space in the case of continuous models without losing the causal information. However, in practice, the geometric model realizes a precubical set. In this case, a natural (but not necessarily optimal) choice is the set of vertices of the precubical set. 

The composition of continuous paths being associative only up to homotopy, we want to use Moore directed paths. Thus from each space of nonconstant directed paths $P^1_{\alpha,\beta}$, we consider the family of reparametrized nonconstant directed paths $\{P^{\ell}_{\alpha,\beta}\mid \ell>0\}$ where $\ell$ is the length. Note that all $P^{\ell}_{\alpha,\beta}$ are homeomorphic to $P_{\alpha,\beta}^1$ for fixed $\alpha,\beta$. Then we consider the family of Moore compositions $P_{\alpha,\beta}^{\ell_1}\p P_{\beta,\gamma}^{\ell_2}\to P_{\alpha,\gamma}^{\ell_1 + \ell_2}$ for all real numbers $\ell_1,\ell_2>0$ and all $\alpha,\beta,\gamma$ belonging to the chosen set of states. It is possible to pack together all these Moore composition maps in an \textit{enriched} semicategorical device which is called a \textit{Moore flow} \cite[Section~6]{Moore1}. The enrichment is necessary to take into account the topology of the space of reparametrization maps. Indeed, the reparametrization must be continuous with respect \textit{both} to the directed path and to the choice of the reparametrization map. One then obtains a \textit{strictly} associative composition law without having to consider directed paths up to nondecreasing reparametrization (these equivalence classes of directed paths are usually called \textit{traces} according to \cite{reparam,reparam-fixed}). 

It is well established that the computer-scientific properties of a concurrent system depend only on the homotopy types of the spaces $P_{\alpha,\beta}^{\ell}$ \cite{DAT_book}. Thanks to their semicategorical nature, Moore flows enable us to prove by pure model categorical arguments (i.e. without explicit calculation) that the space of nonconstant directed paths between two vertices in the geometric realization of a precubical set is m-cofibrant \cite[Corollary~6.8]{RegularMoore}. This fact is originally proved in \cite[Theorem~6.1 and Theorem~7.6]{MR4070250} by constructing an explicit homotopy equivalence with the classifying space of a small category obtained from the precubical set, namely the small category of Ziemia\'{n}ski cube chains associated with the precubical set. As noticed in \cite{Moore1} by one of the anonymous referees, the notion of Moore flow is also an abstraction of the Moore path (semi)category of a topological space, which suggests possible connections with some models of type theory involving Moore paths \cite{MR3903058,north2019typetheoretic}. 

The purpose of this paper is threefold. Firstly, the second paper \cite{Moore2} of this series is technically limited to deal with reparametrization by \textit{nondecreasing homeomorphisms} between nontrivial segments of the real line, instead of with reparametrization by \textit{nondecreasing surjective maps} like in Grandis' notion of $d$-space introduced in \cite{mg}. It is an unnatural restriction which is imposed by the fact that several crucial theorems belonging to the technical core of \cite{Moore2} are either false or their proof is not valid anymore after changing the allowed reparametrizations (cf. Table~\ref{diffGM}). The first purpose of this paper is to fix this issue. Secondly, this paper provides a \textit{uniform treatment} of the two choices of reparametrization setting above by introducing a globular analogue of Raussen's notion of \textit{natural directed path}. It is already known that the cubical version of this notion is central for analyzing the homotopy type of the space of directed paths between two vertices of a precubical set (e.g. \cite[Theorem~6.1 and Theorem~7.6]{MR4070250}). This paper demonstrates that the notion of \textit{natural directed path} is important both for the globular and for the cubical approaches of directed homotopy for concurrency. In fact, we even speculate that this notion is the key for reaching a (still conjectural) unified axiomatic setting which would contain both the globular and cubical approaches of directed homotopy for concurrency. Finally, we want to believe that this work is a contribution in the direction of finding better model categories adapted to directed homotopy for concurrency. The ultimate goal is to find a convenient model category on a category closer to the one of Grandis' $d$-spaces: we would like to remove the multipointed setting somehow. Some speculations about this problem are available in \cite{DHH}. The multipointed setting is a technical restriction introduced in \cite{mdtop} to prevent weak equivalences from contracting the directed segment in the direction of time, because this may destroy the causal structure (as explained above) and therefore this erases the relevant information. There are many speculations about what is a good notion of weak equivalence for a non-multipointed (i.e. continuous) model \cite{zbMATH07374108}. In particular such a notion should be invariant by refinement of observation. There are already techniques to deal with the invariance by refinement of observation in multipointed models, in the cubical setting in \cite{NaturalHomology} and in the globular setting in \cite{3eme,4eme}, which remain to be unified.

The reparametrizations of execution paths allowed in \cite{Moore2} are therefore the precompositions by the maps of the reparametrization category $\mathcal{G}$ in the sense of Proposition~\ref{paramG} which are the nondecreasing homeomorphisms between nontrivial segments of the real line. The technical advantage of this setting is that all execution paths of a \textit{cellular} multipointed $d$-space are \textit{regular} in the sense of \cite[Definition~1.1]{reparam}, namely without \textit{stop intervals}, i.e. without nontrivial intervals on which the path is constant (cf. Definition~\ref{def_regular}). We want to replace the reparametrization category $\mathcal{G}$ by the reparametrization category $\mathcal{M}$ in the sense of Proposition~\ref{paramM} whose maps consist of the nondecreasing surjective maps between nontrivial segments of the real line. The technical obstacle to overcome is that the execution paths of a cellular multipointed $d$-space can now contain stop intervals. By \cite[Proposition~2.2]{reparam-fixed}, every nonconstant \textit{Moore} path in a Hausdorff space has a regular reparametrization. Moreover, by \cite[Proposition~3.8]{reparam}, the regular reparametrizations of a given nonconstant Moore path in a Hausdorff space are unique up to a map of $\mathcal{G}$. In \cite[Definition~2.14]{MR2521708}, Raussen introduces a cubical notion of \textit{naturalization} of a directed path. Intuitively, it means that any nonconstant \textit{directed} path $\gamma$ in the geometric realization of a precubical set has a regular reparametrization called the naturalization $\natcub(\gamma)$ which is morally more natural than the other ones. It is the unique reparametrization which makes the directed path a Moore composition of isometries for some Lawvere metric structure on the geometric realization of a precubical set. This idea is generalized to the setting of symmetric transverse sets in \cite{DirectedDegeneracy} and to the setting of presheaf categories on a thick category of cubes in \cite{ThickCubes}. By Proposition~\ref{local-inj-reg}, there is then a unique factorization $\gamma=\natcub(\gamma)\eta$ where $\eta\in \mathcal{M}$. The technical innovation of this paper is the introduction of a \textit{globular version} of Raussen's idea of naturalization of a directed path in Proposition~\ref{normal-form} and Definition~\ref{def-regularization}. It turns out that every execution path $\gamma$ of a cellular multipointed $d$-space has also a regular reparametrization $\natgl(\gamma)$ which is morally more natural than the other ones. Again by Proposition~\ref{local-inj-reg}, and since the underlying space of a cellular multipointed $d$-space is Hausdorff, there is then a unique factorization $\gamma=\natgl(\gamma)\eta$ where $\eta\in \mathcal{M}$. It is the key point to adapt the technical core of \cite{Moore2}.

Raussen's naturalization and the \textit{globular naturalization} have in common the following property: the naturalization of a Moore composition is the Moore composition of the naturalizations. On the other hand, Raussen's naturalization and the globular naturalization do not behave in the same way with respect to continuous deformations. This point is explained in Corollary~\ref{diff-glob-cube} and in the remark following it. Two directed paths in the geometric realization of a precubical set which are dihomotopy equivalent relatively to the extremities have naturalizations of the same length. On the contrary, the best that can be said in the globular case is that, on such a \textit{compact} continuous path of directed paths, the natural length is bounded (actually, it takes finitely many values). 

The main results of this paper can be stated as follows. The inclusion functor $\mathcal{G}\subset \mathcal{M}$ induces a forgetful functor \[\ptop{\mathcal{M}} \longrightarrow \ptop{\mathcal{G}}\] from $\mathcal{M}$-multipointed $d$-spaces to $\mathcal{G}$-multipointed $d$-spaces and a forgetful functor \[\dtopM \longrightarrow \dtopG\] from $\mathcal{M}$-flows to $\mathcal{G}$-flows. Write $\ptopsat{\mathcal{P}}$ for the category of \textit{saturated} $\mathcal{P}$-multipointed $d$-spaces (the saturation hypothesis is a very important notion of DAT which is introduced in Section~\ref{section-saturated} in the setting of $\mathcal{P}$-multipointed $d$-spaces).

\begin{thmN} (Proposition~\ref{final2}, Theorem~\ref{final}, Theorem~\ref{comparisonGM}, Theorem~\ref{Quillen-equiv-sat} and Corollary~\ref{Quillen-equiv-sat-2})
	There is the commutative diagram of right Quillen equivalences between the five q-model structures
	\[
	\begin{tikzcd}[row sep=5em, column sep=5em]
		\ptopsat{\mathcal{M}} \arrow[r,"\ref{Quillen-equiv-sat}","\subset"'] \arrow[d,"\ref{Quillen-equiv-sat-2}"']& \ptop{\mathcal{M}}\arrow[r,"\ref{comparisonGM}"] \arrow[d,"\ref{final}"'] & \ptop{\mathcal{G}} = \ptopsat{\mathcal{G}} \arrow[d,"\ref{final}"] \\
		\dtopM \arrow[r,equal]& \dtopM  \arrow[r,"\ref{final2}"] & \dtopG.
	\end{tikzcd}
	\]
	Moreover, the unit maps and the counit maps of the three vertical right adjoints induce isomorphisms on q-cofibrant objects.
\end{thmN}

Note that Proposition~\ref{final2} should have been put in \cite{Moore1} as an application of the results of the latter paper: it is an omission. As byproducts of this paper, we also prove the following two results:

\begin{thmN} (Theorem~\ref{main})
	Let $\mathcal{P}$ be either $\mathcal{G}$ or $\mathcal{M}$. The compact-open topology on the set of execution paths of a locally finite cellular $\mathcal{P}$-multipointed $d$-space is $\Delta$-generated. Therefore in this case, the space of execution paths is metrizable with the distance of the uniform convergence.
\end{thmN}

\begin{thmN} (Theorem~\ref{quotient-reparam})
Let $\mathcal{P}$ be either $\mathcal{G}$ or $\mathcal{M}$. Let $X$ be a q-cofibrant $\mathcal{P}$-multipointed $d$-space. Let $\alpha,\beta\in X^0$. Then the quotient map \[\PA_{\alpha,\beta}^{top}X \longrightarrow \PA_{\alpha,\beta} X\] is a homotopy equivalence from an m-cofibrant space to a q-cofibrant space.
\end{thmN}

Theorem~\ref{main} was not in \cite{Moore2}, even for the case $\mathcal{P}=\mathcal{G}$. Theorem~\ref{quotient-reparam} was proved in \cite{model2} for the case $\mathcal{P} = \mathcal{G}$, but the proof does not seem to be generalizable to the case $\mathcal{P}=\mathcal{M}$ (see the long comment before the statement of the theorem).

In addition to generalizing the results of \cite{Moore2} and to finding proofs which are independent of the choice of the reparametrization category $\mathcal{G}$ or $\mathcal{M}$ (except for Theorem~\ref{img-closed}: the statement is independent of the choice of $\mathcal{G}$ or $\mathcal{M}$, but not the proof), this work raises the question of finding a better definition of a reparametrization category than Definition~\ref{def-reparam}. We suspect that there is a model category containing the reparametrization categories such that $\mathcal{G}$ and $\mathcal{M}$ are cofibrant replacements of the terminal category. A cofibrant replacement $\mathcal{P}$ of the terminal category in this hypothetical model category should give rise to a notion of $\mathcal{P}$-multipointed $d$-space.

\subsection*{Outline of the paper}

Section~\ref{Mooreflow} is a reminder about $\mathcal{P}$-flows for a reparametrization category $\mathcal{P}$ which is either $\mathcal{G}$ or $\mathcal{M}$ in this paper. Section~\ref{multipointedspace} adapts some results and constructions for $\mathcal{G}$-multipointed $d$-spaces proved in \cite{Moore2} to the case of $\mathcal{P}$-multipointed $d$-spaces. Section~\ref{core} is the adaptation of \cite[Section~5]{Moore2} to the case of cellular $\mathcal{P}$-multipointed $d$-spaces. The main results are the notion of globular naturalization of an execution path of a cellular $\mathcal{P}$-multipointed $d$-space (Proposition~\ref{normal-form} and Definition~\ref{def-regularization}). We then obtain, thanks to the notion of carrier of an execution path, Theorem~\ref{calcul_final_structure} which is a replacement for \cite[Theorem~5.20]{Moore2} and Theorem~\ref{carrier-finite-on-compact} which is a replacement for \cite[Theorem~5.19]{Moore2}. Section~\ref{core2} is the adaptation of \cite[Section~6]{Moore2}. A generalization of \cite[Theorem~6.11]{Moore2} is proved in Theorem~\ref{img-closed}. Section~\ref{locallyfinite} is a digression which uses Theorem~\ref{diagonal-execution} and Theorem~\ref{img-closed} to prove that the space of execution paths in the locally finite case is metrizable with the distance of the uniform convergence in Theorem~\ref{main}. Section~\ref{unit} establishes the main theorems of the paper, namely Theorem~\ref{final} and Theorem~\ref{comparisonGM}. Section~\ref{section-saturated} proves that the saturation hypothesis (it is meaningless for the case $\mathcal{P}=\mathcal{G}$ by Proposition~\ref{ex-non-saturated}), which is a very important notion in DAT, can be safely added to the definition of an $\mathcal{M}$-multipointed $d$-space without changing the mathematical properties.

\subsection*{Erratum}

As explained in the corrected version of \cite{Moore1}, the tenseur product of $\mathcal{P}$-spaces is not symmetric if $\mathcal{P}$ is $\mathcal{G}$ or $\mathcal{M}$. Therefore, the word symmetric must be removed everywhere from \cite{Moore2}. Besides, the terminology of \textit{biclosed} semimonoidal structure should be used instead of the terminology of closed semimonoidal structure to describe the tensor product of $\mathcal{P}$-spaces.

\subsection*{Prerequisites and notations}

We refer to \cite{TheBook} for locally presentable categories, to \cite{MR2506258} for combinatorial model categories.  We refer to \cite{MR99h:55031,ref_model2} for more general model categories. We refer to \cite{KellyEnriched} and to \cite[Chapter~6]{Borceux2} for enriched categories. All enriched categories are topologically enriched categories: \textit{the word topologically is therefore omitted}. A gold mine of examples and counterexamples in general topology can be found in \cite{zbMATH06070728}. A \textit{cellular object} of a combinatorial model category is an object $X$ such that the canonical map $\varnothing\to X$ is a transfinite composition of pushouts of generating cofibrations.

The results of this paper rely heavily on the results of \cite{Moore2}. A self-contained paper would not help the reader much. The choice made for this work is to emphasize the differences between $\mathcal{G}$ and $\mathcal{M}$ instead of the similarities. Table~\ref{diffGM} summarizes these differences. The left column is a list of theorems of \cite{Moore2}. The middle column gives the status of the statement for $\mathcal{P}=\mathcal{M}$. The right column gives the replacement in this paper: it consists of a statement which is modified if necessary and a new proof. In this paper, even if \cite[Theorem~5.19]{Moore2} is still valid for $\mathcal{P}=\mathcal{M}$, it is replaced by Theorem~\ref{carrier-finite-on-compact} which is a much powerful statement both for the proof of Theorem~\ref{pre-calculation-pathspace} and to understand the difference between the globular naturalization and the cubical naturalization.

The category $\topspace$ denotes the category of \textit{$\Delta$-generated spaces} or of \textit{$\Delta$-Hausdorff $\Delta$-generated spaces} (cf. \cite[Section~2 and Appendix~B]{leftproperflow}). The inclusion functor from the full subcategory of $\Delta$-generated spaces to the category of general topological spaces together with the continuous maps has a right adjoint called the $\Delta$-kelleyfication functor. The latter functor does not change the underlying set: it only adds open subsets. The category $\topspace$ is locally presentable (see \cite[Corollary~3.7]{FR} for the non $\Delta$-Hausdorff case and \cite[Proposition~B.18]{leftproperflow} for the $\Delta$-Hausdorff case),  and cartesian closed by a theorem of Dugger-Vogt recalled in \cite[Proposition~2.5]{mdtop}. The internal hom $\ttop(X,Y)$ is given by taking the $\Delta$-kelleyfication of the compact-open topology on the set $\topspace(X,Y)$. The category $\topspace$ is equipped with its q-model structure denoted by $\topspace_q$. The h-model structure of $\topspace$ provided by \cite[Corollary~5.23]{Barthel-Riel} is mentioned before Theorem~\ref{quotient-reparam}. The m-model structure of $\topspace$ in the sense of \cite{mixed-cole} is briefly used in the proof of Theorem~\ref{quotient-reparam}. \textit{A compact space is a quasicompact Hausdorff space (French convention)}. All $\Delta$-generated spaces are sequential. 

The proofs of this paper rely on the facts that the $\Delta$-generated spaces are colimits of the segment $[0,1]$ and that all involved topological spaces are sequential. Some intermediate steps even use sequential spaces which are not necessarily $\Delta$-generated:  in Theorem~\ref{pre-calculation-pathspace}, we need to construct a finite covering of $[0,1]$ by closed subsets equipped with the relative topology. However, using the same techniques as in \cite[Appendix~C]{Moore2}, the main theorems of this paper (the Quillen equivalences and Theorem~\ref{GM-cat}) can be extended to other convenient categories of topological spaces for doing algebraic topology like $k$-spaces.

$\K^{op}$ denotes the opposite category of $\K$; $\Obj(\K)$ is the class of objects of $\K$; $\K^I$ is the category of functors and natural transformations from a small category $I$ to $\K$; $\varnothing$ is the initial object, $\mathbf{1}$ is the final object, $\id_X$ is the identity of $X$; $\K(X,Y)$ is the set of maps in a set-enriched, i.e. locally small, category $\K$; $\K(X,Y)$ is the space of maps in an enriched category $\K$. The underlying set of maps may be denoted by $\K_0(X,Y)$ if it is necessary to specify that we are considering the underlying set. Table~\ref{CheatSheet} is a general overview for some other notations used in this paper.

\textit{All Moore paths in this paper are nonconstant}: see Definition~\ref{def_Moore_path}.

\begin{table}
	\begin{tabular}{|c||c|c|}
		\hline
		 & Status for $\mathcal{M}$ & Replacement \\
		\hline\hline
		\cite[Theorem~3.9]{Moore2} & true & Theorem~\ref{final-structure-revisited} \\
		\hline
		\cite[Proposition~2.12]{Moore2} & true & Proposition~\ref{calcul-topology-glob} \\
		\hline
		\cite[Theorem~5.7]{Moore2} & true & Theorem~\ref{cof-accessible} \\
		\hline
		\cite[Theorem~5.9]{Moore2} & false (1) & Theorem~\ref{normal-form} \\
		\hline
		\cite[Theorem~5.20]{Moore2} & false (2) & Theorem~\ref{calcul_final_structure} \\
		\hline
		\cite[Proposition~5.17]{Moore2} & false (3) & Theorem~\ref{diagonal-execution} \\
		\hline
		\cite[Theorem~5.18]{Moore2} & true & Corollary~\ref{bounded0} \\
		\hline
		\cite[Theorem~5.19]{Moore2} & true & Corollary~\ref{bounded} \\
		\hline
		\cite[Proposition~6.3]{Moore2} & true & Proposition~\ref{comp-gl} \\
		\hline
		\cite[Theorem~6.11]{Moore2} & true & Theorem~\ref{img-closed} \\
		\hline
		\cite[Theorem~7.2 and Theorem~7.3]{Moore2} & true & Theorem~\ref{pre-calculation-pathspace}\\
		\hline
	\end{tabular}
	\caption{Main differences between $\mathcal{G}$ and $\mathcal{M}$: (1) see the comment after Theorem~\ref{normal-form}, (2) see the comment after Theorem~\ref{calcul_final_structure}, (3) see Proposition~\ref{impossible}}
	\label{diffGM}
\end{table}

\begin{table}
	\begin{tabular}{|c|c|}
		\hline
		Space of \textit{all} execution paths of  &The space $\PA^{top}X$\\length $1$ of a $\mathcal{P}$-multipointed $d$-space $X$ & \\
		\hdashline
		And only from $\alpha$ to $\beta$ & The space $\PA^{top}_{\alpha,\beta}X$\\
		&(the superscript $top$ means that the \\&execution  paths  are continuous paths)\\
		\hline
		Space of \textit{all} execution paths of  &The space $\PA^{\ell}X$\\ length $\ell$ of a $\mathcal{P}$-multipointed $d$-space $X$ & \\
		\hdashline
		And only from $\alpha$ to $\beta$ & The space $\PA^{\ell}_{\alpha,\beta}X$\\
		&(the superscript $top$ is replaced \\& by the length)\\
		\hline
		$\mathcal{P}$-space of \textit{all} execution paths &  The $\mathcal{P}$-space $\PA X$\\of a $\mathcal{P}$-flow $X$ &\\
		\hdashline
		And only from $\alpha$ to $\beta$ & The $\mathcal{P}$-space $\PA_{\alpha,\beta}X$\\
		&(the reparametrization category \\&is determined by $X$)\\
		\hline
		Space of \textit{all} execution paths& The space $\PA X$\\ of a flow $X$ & \\
		\hdashline
		And only from $\alpha$ to $\beta$ & The space $\PA_{\alpha,\beta}X$\\
		&(the reparametrization category \\& is here the terminal category)\\
		\hline
		Topological globe of a space $Z$ & The $\mathcal{P}$-multipointed $d$-space $\globP(Z)$\\ 
		& (the use of the superscript $\mathcal{P}$  is necessary \\&   to specify the reparametrization category, \\& the information being not in $Z$)\\
		\hline
		Globe of a $\mathcal{P}$-space $Z$ & The $\mathcal{P}$-flow $\glob(Z)$ \\ 
		& (the reparametrization category \\&is determined by $Z$ which is a $\mathcal{P}$-space)\\
		\hline
		Globe of a space $Z$ & The flow $\glob(Z)$ \\& (the reparametrization category \\& is here the terminal category)\\
		\hline
	\end{tabular}
	\caption{Overview of some notations used in this paper: see also Remark~\ref{rem:cheatsheet}}
	\label{CheatSheet}
\end{table}

\subsection*{Acknowledgments}

I would like to thank the anonymous referee for many helpful comments, which have improved the presentation of the paper. I would also like to thank the editors for their help in the final preparation of this document.

\section{Moore flow}
\label{Mooreflow}

\begin{nota}
		The notations $\ell,\ell',\ell_i,L,\dots$ mean a strictly positive real number unless specified something else. $[\ell,\ell']$ denotes a segment: unless specified, it is always understood that $\ell<\ell'$.
\end{nota}

\bd \label{def-reparam} \cite[Definition~4.3]{Moore1}
A \textit{reparametrization category} $(\mathcal{P},\ot)$ is a small enriched semimonoidal category satisfying the following additional properties: 
\begin{enumerate}
	\item The semimonoidal structure is strict, i.e. the associator is the identity.
	\item All spaces of maps $\mathcal{P}(\ell,\ell')$ for all objects $\ell$ and $\ell'$ of $\mathcal{P}$ are contractible. 
	\item For all maps $\phi:\ell\to \ell'$ of $\mathcal{P}$, for all $\ell'_1,\ell'_2\in \Obj(\mathcal{P})$ such that $\ell'_1\ot\ell'_2=\ell'$, there exist two maps $\phi_1:\ell_1\to \ell'_1$ and $\phi_2:\ell_2\to \ell'_2$ of $\mathcal{P}$ such that $\phi=\phi_1 \ot \phi_2 : \ell_1\ot\ell_2 \to \ell'_1 \ot\ell'_2$ (which implies that $\ell_1 \ot \ell_2=\ell$). 
\end{enumerate}
\ed 

The terminal category is a symmetric reparametrization category. It is not known whether there exist symmetric reparametrization categories not equivalent to the terminal category. Here are the two examples of reparametrization category used in this paper.

\bp \label{paramG} \cite[Proposition~4.9]{Moore1}
There exists a reparametrization category, denoted by ${\mathcal{G}}$, such that the semigroup of objects is the open interval $]0,+\infty[$ equipped with the addition and such that for every $\ell_1,\ell_2>0$, $\mathcal{G}(\ell_1,\ell_2)$ is the set of nondecreasing homeomorphisms from $[0,\ell_1]$ to $[0,\ell_2]$ equipped with the $\Delta$-kelleyfication of the relative topology induced by the set inclusion $\mathcal{G}(\ell_1,\ell_2) \subset \ttop([0,\ell_1],[0,\ell_2])$ and such that for every $\ell_1,\ell_2,\ell_3>0$, the composition map $\mathcal{G}(\ell_1,\ell_2)\p \mathcal{G}(\ell_2,\ell_3) \to \mathcal{G}(\ell_1,\ell_3)$ is induced by the composition of continuous maps.
\ep

\bp \label{paramM} \cite[Proposition~4.11]{Moore1}
There exists a reparametrization category, denoted by ${\mathcal{M}}$, such that the semigroup of objects is the open interval $]0,+\infty[$ equipped with the addition and such that for every $\ell_1,\ell_2>0$, $\mathcal{M}(\ell_1,\ell_2)$ is the set of nondecreasing surjective maps from $[0,\ell_1]$ to $[0,\ell_2]$ equipped with the $\Delta$-kelleyfication of the relative topology induced by the set inclusion $\mathcal{M}(\ell_1,\ell_2) \subset \ttop([0,\ell_1],[0,\ell_2])$ and such that for every $\ell_1,\ell_2,\ell_3>0$, the composition map $\mathcal{M}(\ell_1,\ell_2)\p \mathcal{M}(\ell_2,\ell_3) \to \mathcal{M}(\ell_1,\ell_3)$ is induced by the composition of continuous maps.
\ep

\begin{nota}
	A reparametrization category $\mathcal{P}$ which is either $\mathcal{G}$ or $\mathcal{M}$ is fixed for the rest of the paper. 
\end{nota}

\bp \label{morphG-metrizable}
The topology of $\mathcal{P}(\ell_1,\ell_2)$ is the compact-open topology. In particular, it is metrizable. A sequence $(\phi_n)_{n\geq 0}$ of $\mathcal{P}(\ell_1,\ell_2)$ converges to $\phi\in \mathcal{P}(\ell_1,\ell_2)$ if and only if it converges pointwise. 
\ep

\bpf
It is mutatis mutandis the same argument as the one given for $\mathcal{P}=\mathcal{G}$ in \cite[Proposition~2.5]{Moore2}. 
\epf

\begin{nota}
	Let $\phi_i\in \mathcal{P}(\ell_i,\ell'_i)$ for $n\geq 1$ and $1\leqslant i \leqslant n$. Then the map
	\[
	\phi_1 \ot \dots \ot \phi_n : \sum_i \ell_i \longrightarrow \sum_i \ell'_i
	\]
	denotes the nondecreasing surjective map defined by 
	\[
	(\phi_1 \ot \dots \ot \phi_n)(t) = \begin{cases}
		\phi_1(t) & \hbox{if } 0\leqslant t\leqslant \ell_1\\
		\phi_2(t-\ell_1)+\ell'_1 & \hbox{if } \ell_1\leqslant t\leqslant \ell_1+\ell_2\\
		\dots \\
		\phi_i(t-\sum_{j<i}\ell_j) + \sum_{j<i}\ell'_j& \hbox{if } \sum_{j<i}\ell_j\leqslant t \leqslant \sum_{j\leqslant i}\ell_j\\
		\dots\\
		\phi_n(t-\sum_{j<n}\ell_j) + \sum_{j<n}\ell'_j & \hbox{if } \sum_{j<n}\ell_j\leqslant t \leqslant \sum_{j\leqslant n}\ell_j.
	\end{cases}
	\] 
\end{nota}

\bp \label{decomposition-tenseur}
Let $\phi\in \mathcal{P}(\ell,\ell')$. Let $n\geq 1$. Consider $\ell'_1,\dots,\ell'_n>0$ with $n\geq 1$ such that $\ell'_1 + \dots + \ell'_n = \ell'$. Then there exists a decomposition of $\phi$ of the form \[\phi=\phi_1\ot \dots \ot \phi_n\] such that $\phi_i\in \mathcal{P}(\ell_i,\ell'_i)$ for $1\leqslant i \leqslant n$. Moreover, if $\mathcal{P}=\mathcal{G}$, then this decomposition is unique.
\ep

\bpf
The case $n=1$ is trivial. The case $n=2$ comes from the fact $\mathcal{G}$ and $\mathcal{M}$ are reparametrization categories. We deduce the existence of the decomposition by induction on $n\geq 2$. The uniqueness when $\mathcal{P}=\mathcal{G}$ is \cite[Proposition~3.2]{Moore2}.
\epf

\begin{nota} \label{Rspace}
	The enriched category of enriched presheaves from $\mathcal{P}$ to $\topspace$ is denoted by $\topdgrP$. The underlying set-enriched category of enriched maps of enriched presheaves is denoted by $\topdgrP_0$. The objects of $\topdgrP_0$ are called the \textit{$\mathcal{P}$-spaces}. Let \[\mathbb{F}^{\mathcal{P}^{op}}_{\ell}U=\mathcal{P}(-,\ell)\p U \in \topdgrP_0\] where $U$ is a topological space and where $\ell>0$.
\end{nota}

\bp\label{ev-adj} \cite[Proposition~5.3 and Proposition~5.5]{dgrtop}
The category $\topdgrP_0$ is a full reflective and coreflective subcategory of $\topspace^{\mathcal{P}^{op}_0}$. 
For every $\mathcal{P}$-space $F:\mathcal{P}^{op}\to \topspace$, every $\ell>0$ and every topological space $X$, we have the natural bijection of sets \[\topdgrP_0(\mathbb{F}^{\mathcal{P}^{op}}_{\ell}X,F) \iso \topspace(X,F(\ell)).\] 
\ep

\bth (\cite[Proposition~5.1]{dgrtop} and \cite[Theorem~5.14]{Moore1}) \label{closedsemimonoidal}
The category $\topdgrP_0$ is locally presentable. Let $D$ and $E$ be two $\mathcal{P}$-spaces. Let 
\[
D \ot E = \int^{(\ell_1,\ell_2)} \mathcal{P}(-,\ell_1+\ell_2) \p D(\ell_1) \p E(\ell_2).
\]
The pair $(\topdgrP_0,\ot)$ has the structure of a biclosed semimonoidal category.
\eth

\bd \cite[Definition~6.2]{Moore1} \label{def-Moore-flow}
A \textit{$\mathcal{P}$-flow}, also called a \textit{Moore flow} if there is no ambiguity on the choice of $\mathcal{P}$, is a small semicategory enriched over the biclosed semimonoidal category $(\topdgrP_0,\ot)$ of Theorem~\ref{closedsemimonoidal}. The corresponding category is denoted by $\dtopP$. 
\ed

A $\mathcal{P}$-flow $X$ consists of a \textit{set of states} $X^0$, for each pair $(\alpha,\beta)$ of states a $\mathcal{P}$-space $\PA_{\alpha,\beta}X$ of $\topdgrP_0$ and for each triple $(\alpha,\beta,\gamma)$ of states an associative composition law \[*:\PA_{\alpha,\beta}X \ot \PA_{\beta,\gamma}X \to \PA_{\alpha,\gamma}X.\] A map of $\mathcal{P}$-flows $f$ from $X$ to $Y$ consists of a set map \[f^0:X^0 \to Y^0\] (often denoted by $f$ as well if there is no possible confusion) together for each pair of states $(\alpha,\beta)$ of $X$ with a natural transformation \[\PA f:\PA_{\alpha,\beta}X \longrightarrow \PA_{f(\alpha),f(\beta)}Y\] compatible with the composition law. The topological space $\PA_{\alpha,\beta}X(\ell)$ is denoted by $\PA_{\alpha,\beta}^\ell X$ and is called the space of \textit{execution paths of length $\ell$}.

\bd \label{Pmoore}
	Let $X$ be a $\mathcal{P}$-flow. The \textit{$\mathcal{P}$-space of execution paths} $\PA X$ of $X$ is by definition the $\mathcal{P}$-space \[\PA X = \coprod_{(\alpha,\beta)\in X^0\p X^0} \PA_{\alpha,\beta}X.\]
	It yields a well-defined functor $\PA:\dtopP \to \topdgrP_0$. The image of $\ell$ is denoted by $\PA^\ell$. We therefore have the equality \[\PA^\ell X = \coprod_{(\alpha,\beta)\in X^0\p X^0} \PA_{\alpha,\beta}^{\ell}X.\]
\ed

The category $\dtopP$ is locally presentable by \cite[Theorem~6.11]{Moore1}. 

\begin{nota}
	Let $D:\mathcal{P}^{op}\to \topspace$ be a $\mathcal{P}$-space. We denote by $\glob(D)$ the Moore flow defined as follows: 
	\[
	\begin{aligned}
		&\glob(D)^0 = \{0,1\}\\
		&\PA_{0,0}\glob(D)=\PA_{1,1}\glob(D)=\PA_{1,0}\glob(D)=\varnothing\\
		&\PA_{0,1}\glob(D)=D.
	\end{aligned}	
	\]
	There is no composition law. This construction yields a functor \[\glob:\topdgrP_0\to \dtopP.\]
\end{nota}

\begin{rem}
	The notation $\glob(D)$ is not ambiguous since $D$ is always a $\mathcal{P}$-space with $\mathcal{P}$ being $\mathcal{G}$ or $\mathcal{M}$ and $\glob(D)$ is then necessarily either a $\mathcal{G}$-flow or an $\mathcal{M}$-flow respectively.
\end{rem}

By \cite[Theorem~6.2]{dgrtop}, the category of $\mathcal{P}$-spaces $\topdgrP_0$ can be endowed with the projective model structure associated with the model structure $\topspace_q$. It is called the projective q-model structure. It is combinatorial. The fibrations are the objectwise q-fibrations. The weak equivalences are the objectwise weak homotopy equivalences. All $\mathcal{P}$-spaces are fibrant for this model structure. By \cite[Theorem~8.8, Theorem~8.9 and Theorem~8.16]{Moore1}, the category of $\mathcal{P}$-flows can be endowed with a combinatorial model structure characterized as follows: 1) a map of $\mathcal{P}$-flows $f:X\to Y$ is a weak equivalence if and only if $f^0:X^0\to Y^0$ is a bijection and $\PA f:\PA_{\alpha,\beta}X\to \PA_{f(\alpha),f(\beta)}Y$ is a weak equivalence of the projective q-model structure of $\mathcal{P}$-flows; 2) a map of $\mathcal{P}$-flows $f:X\to Y$ is a fibration if and only if $\PA f:\PA_{\alpha,\beta}X\to \PA_{f(\alpha),f(\beta)}Y$ is a fibration of the projective q-model structure of $\mathcal{P}$-flows, i.e. an objectwise q-fibration of topological spaces. All $\mathcal{P}$-flows are q-fibrant.

\section{Multipointed d-space}
\label{multipointedspace}

\bd \label{composition_map} Let $\gamma_1$ and $\gamma_2$ be two continuous maps from $[0,1]$ to some topological space such that $\gamma_1(1)=\gamma_2(0)$. The composite defined by 
\[
(\gamma_1 *_N \gamma_2)(t) = \begin{cases}
	\gamma_1(2t)& \hbox{ if }0\leqslant t\leqslant \frac{1}{2},\\
	\gamma_2(2t-1)& \hbox{ if }\frac{1}{2}\leqslant t\leqslant 1
\end{cases}
\]
is called the \textit{normalized composition}. The normalized composition being not associative, a notation like $\gamma_1 *_N \dots *_N \gamma_n$ will mean, by convention, that $*_N$ is applied from the left to the right. 
\ed

\bd \label{def_Moore_path}
Let $U$ be a topological space. A \textit{(Moore) path} of $U$ consists in this paper of a \textit{nonconstant} continuous map $\gamma:[0,\ell]\to U$ with $\ell>0$. The real number $\ell$ is called the \textit{length} of $\gamma$.
\ed

Let $\gamma_1:[0,\ell_1]\to U$ and $\gamma_2:[0,\ell_2]\to U$ be two Moore paths of a topological space $U$ such that $\gamma_1(\ell_1)=\gamma_2(0)$. The \textit{Moore composition} $\gamma_1*\gamma_2:[0,\ell_1+\ell_2]\to U$ is the Moore path defined by 
\[
(\gamma_1*\gamma_2)(t)=
\begin{cases}
	\gamma_1(t) & \hbox{ for } t\in [0,\ell_1]\\
	\gamma_2(t-\ell_1) &\hbox{ for }t\in [\ell_1,\ell_1+\ell_2].
\end{cases}
\]
The Moore composition of Moore paths is strictly associative.

\begin{nota}
	Let $\ell>0$. Let $\mu_{\ell}:[0,\ell]\to [0,1]$ be the homeomorphism defined by $\mu_\ell(t) = t/\ell$.
\end{nota}

\bd A \textit{$\mathcal{P}$-multipointed $d$-space $X$} or just \textit{multipointed $d$-space $X$} if there is no ambiguity on the choice of $\mathcal{P}$ is a triple $(|X|,X^0,\PA^{top}X)$ where
\begin{itemize}[leftmargin=*]
	\item The pair $(|X|,X^0)$ is a multipointed space. The space $|X|$ is called the \textit{underlying space} of $X$ and the set $X^0$ the \textit{set of states} of $X$.
	\item The set $\PA^{top}X$ is a set of continuous maps from $[0,1]$ to $|X|$ called the \textit{execution paths}, satisfying the following axioms:
	\begin{itemize}
		\item For any execution path $\gamma$, one has $\gamma(0),\gamma(1)\in X^0$.
		\item Let $\gamma$ be an execution path of $X$. Then any composite $\gamma\phi$ with $\phi\in \mathcal{P}(1,1)$ is an execution path of $X$.
		\item Let $\gamma_1$ and $\gamma_2$ be two composable execution paths of $X$; then the normalized composition $\gamma_1 *_N \gamma_2$ is an execution path of $X$.
	\end{itemize}
\end{itemize}
A map $f:X\to Y$ of $\mathcal{P}$-multipointed $d$-spaces is a map of multipointed spaces from $(|X|,X^0)$ to $(|Y|,Y^0)$ such that for any execution path $\gamma$ of $X$, the map $\PA^{top}f:\gamma\mapsto f. \gamma$ is an execution path of $Y$. 
\ed

\begin{nota}
	The category of $\mathcal{P}$-multipointed $d$-spaces is denoted by $\ptop{\mathcal{P}}$. 
\end{nota}

Proposition~\ref{locally-presentable} is already known for $\mathcal{P} = \mathcal{G}$ by \cite[Theorem~3.5]{mdtop} when $\topspace$ is the category of $\Delta$-generated spaces. We recall some details and we give some references to make the argument of Proposition~\ref{saturated-locally-presentable} understandable for a reader not familiar with the logical characterization of locally presentable categories.

\bp \label{locally-presentable}
The category $\ptop{\mathcal{P}}$ is locally presentable.
\ep

\bpf
We use the terminology of \cite[Chapter~5]{TheBook}. Let $\mathcal{T}$ be a small relational universal strict Horn theory axiomatizing $\topspace$ (without equality by \cite[Theorem~3.6]{FR} when $\topspace$ is the category of $\Delta$-generated spaces, and with equality by \cite[Proposition~B.18]{leftproperflow} when $\topspace$ is the category of $\Delta$-Hausdorff $\Delta$-generated spaces). Let $\{R_j\mid j\in J\}$ be the set of relational symbols of $\mathcal{T}$. We add now to $\mathcal{T}$ the axioms encoding the structure of $\mathcal{P}$-multipointed $d$-space on a topological space. We need a $1$-ary relational symbol $S$ whose interpretation is the set of states and a $[0,1]$-ary relational symbol $R$ whose interpretation is the set of execution paths. The axioms look like as follows (the first axiom appears in the proof of \cite[Theorem~3.5]{mdtop} and the other axioms are already in the proof of \cite[Theorem~4.2]{FR}):
\begin{enumerate}
	\item $(\forall x) R(x) \Rightarrow (S(x_0) \wedge S(x_1))$
	\item $(\forall x,y,z) \left(\left(\bigwedge_{0\leqslant t\leqslant 1/2} x_{2t}=z_t\right) \wedge
	\left(\bigwedge_{1/2\leqslant t\leqslant 1} y_{2t-1}=z_{\frac{1}{2}+t}\right) \wedge R(x)
	\wedge R(y)\right) \Rightarrow R(z)$
	\item $(\forall x) R(x) \Rightarrow R(x.t)$ where $t\in \mathcal{P}(1,1)$
	\item $(\forall x) R(x) \Rightarrow R_j(x.a)$ where $j\in J$ and
	$\mathcal{T}$ satisfies $R_j$ for a sequence $a$ of $[0,1]$.
\end{enumerate}
The first three statements encode the structure of $\mathcal{P}$-multipointed $d$-space whereas the fourth one encodes the continuity of execution paths: for example, if the isomorphism $\topspace\iso \mathrm{Mod}(\mathcal{T})$ between $\topspace$ and the models of $\mathcal{T}$ takes the segment $[0,1]$ to the set $[0,1]$ such that $R_{j_0}(1/3,0.2,1/\pi)$ holds for some $j_0\in J$, then $(\forall x) R(x)\Rightarrow R_{j_0}(x_{1/3},x_{0.2},x_{1/\pi})$ is an axiom. By \cite[Theorem~5.30]{TheBook}, the proof is complete. 
\epf

The subset of execution paths from $\alpha$ to $\beta$ is the set of $\gamma\in\PA^{top} X$  such that $\gamma(0)=\alpha$ and $\gamma(1)=\beta$; it is denoted by $\PA^{top}_{\alpha,\beta} X$: $\alpha$ is called the \textit{initial state} and $\beta$ the \textit{final state} of such a $\gamma$. The set $\PA^{top}_{\alpha,\beta} X$ is equipped with the $\Delta$-kelleyfication of the relative topology induced by the inclusion $\PA^{top}_{\alpha,\beta} X\subset \ttop([0,1],|X|)$. It induces a functor \[\PA^{top}:\ptop{\mathcal{P}} \longrightarrow \topspace.\] Unless specified, the set $\PA^{top}_{\alpha,\beta} X$ is always equipped with this topology.

\begin{nota}
	The mapping $\PA^{top}f$ will be often denoted by $f$ if there is no ambiguity.
\end{nota}

The following examples play an important role in the sequel. 
\begin{enumerate}[leftmargin=*]
	\item Any set $E$ will be identified with the $\mathcal{P}$-multipointed $d$-space $(E,E,\varnothing)$.
	\item The \textit{topological globe of $Z$}, which is denoted by $\globP(Z)$, is the $\mathcal{P}$-multipointed $d$-space defined as follows
	\begin{itemize}
		\item the underlying topological space is the quotient space \[\frac{\{{0},{1}\}\sqcup (Z\p[0,1])}{(z,0)=(z',0)={0},(z,1)=(z',1)={1}}\]
		\item the set of states is $\{{0},{1}\}$ ($0$ is called the \textit{initial state} and $1$ the \textit{final state} of the globe)
		\item the set of execution paths is the set of continuous maps \[\{\delta_z\phi\mid \phi\in \mathcal{P}(1,1),z\in  Z\}\]
		with $\delta_z(t) = (z,t)$.	It is equal to the underlying set of the space $\mathcal{P}(1,1)\p Z$.
	\end{itemize}
	In particular, $\globP(\varnothing)$ is the $\mathcal{P}$-multipointed $d$-space $\{{0},{1}\} = (\{{0},{1}\},\{{0},{1}\},\varnothing)$. 
	\item The \textit{directed segment} is the $\mathcal{P}$-multipointed $d$-space $\vI^{\mathcal{P}}=\globP(\{0\})$. 
\end{enumerate}

\begin{rem} \label{rem:cheatsheet} (see also Table~\ref{CheatSheet})
	The terminology \textit{$\mathcal{P}$-multipointed $d$-space} is chosen because it is a variant of the notion of Grandis' $d$-space. The adjective \textit{$\mathcal{P}$-multipointed} contains the information about the allowed reparametrizations and the constraints on the extremities for the execution paths. The set of execution paths of length $1$ is denoted by $\PA^{top}X$ for a $\mathcal{P}$-multipointed $d$-space $X$, and not $\PA^{\mathcal{P}}X$, because it is always a set of continuous paths equipped with the $\Delta$-kelleyfication of the relative topology, and \textit{never} a $\mathcal{P}$-space: it is the meaning of the superscript $top$. On the contrary, for $\mathcal{P}$-flows (see Definition~\ref{def-Moore-flow}) of for flows (see Definition~\ref{def:flow}), the execution paths are not continuous paths in a topological space. Notations like $\PA X$ or $\PA_{\alpha,\beta}X$ are used for the $\mathcal{P}$-space of execution paths when $X$ is a $\mathcal{P}$-flow for a reparametrization category $\mathcal{P}$. The notation $\globP(Z)$ for a given topological space $Z$ is used to specify that $\globP(Z)$ is a $\mathcal{P}$-multipointed $d$-space. The $\mathcal{G}$-multipointed $d$-space $\globG(Z)$ is not equal to the $\mathcal{M}$-multipointed $d$-space $\globM(Z)$ indeed because $\globG(Z)$ and $\globM(Z)$ do not contain the same execution paths. In particular, all execution paths of $\globG(Z)$ are regular.
\end{rem}

\bp
Let $X$ be a $\mathcal{P}$-multipointed $d$-space. Let $(\alpha,\beta)\in X^0\p X^0$. The following data assemble into a $\mathcal{P}$-space denoted by $\PA_{\alpha,\beta}^\bullet X$:
\begin{itemize}
	\item $\PA_{\alpha,\beta}^\ell X = \{\gamma\mu_\ell \mid \gamma\in \PA_{\alpha,\beta}^{top}X\}$
	\item For $\phi:\ell'\to \ell \in \mathcal{P}$ and $\gamma\in \PA_{\alpha,\beta}^\ell X$, $\PA_{\alpha,\beta}^\phi X(\gamma)=\gamma\phi$.
\end{itemize}
\ep

\bpf
Let $\gamma \in \PA_{\alpha,\beta}^\ell X$. Then, by definition of $\PA_{\alpha,\beta}^\ell X$, there exists a (unique) $\overline{\gamma}\in \PA_{\alpha,\beta}^{top}X$ such that $\gamma=\overline{\gamma}\mu_\ell$. We obtain $\gamma\phi=\overline{\gamma}(\mu_\ell \phi \mu_{\ell'}^{-1}) \mu_{\ell'}$. Since $\mu_\ell \phi \mu_{\ell'}^{-1}\in \mathcal{P}(1,1)$, we have $\overline{\gamma}(\mu_\ell \phi \mu_{\ell'}^{-1})\in \PA_{\alpha,\beta}^{top}X$ and therefore $\gamma\phi\in \PA_{\alpha,\beta}^{\ell'} X$. 
\epf

\begin{nota} \label{trace-space}
	Let $X$ be a $\mathcal{P}$-multipointed $d$-space. Let $(\alpha,\beta)\in X^0\p X^0$. Let \[\PA_{\alpha,\beta} X = \liminj \PA_{\alpha,\beta}^\bullet X.\]
\end{nota}

The topological space $\PA_{\alpha,\beta} X = \liminj \PA_{\alpha,\beta}^\bullet X$ is the quotient of the topological space \[\coprod_{\ell>0}\PA_{\alpha,\beta}^\ell X\] by the equivalence relation generated by the identifications $\gamma \sim\gamma' \Leftrightarrow \gamma\phi = \gamma\phi'$ with $\gamma\in \PA_{\alpha,\beta}^\ell X$, $\gamma'\in \PA_{\alpha,\beta}^{\ell'} X$, $\phi\in \mathcal{P}(\ell'',\ell)$ and $\phi'\in \mathcal{P}(\ell'',\ell')$.

\bd (\cite[Definition~1.2]{reparam}) \label{reparam-equivalence} The above equivalence relation $\sim$ is called the \textit{reparametrization equivalence}. The two Moore paths $\gamma$ and $\gamma'$ above are said \textit{reparametri\-zation equivalent}.
\ed

\bp \label{pre-un}
For every $u,v\in \mathcal{P}(1,\ell)$, there exist $\phi_u,\phi_v\in \mathcal{P}(1,1)$ such that $u\phi_u = v\phi_v$.
\ep

\bpf
Let $u,v\in \mathcal{P}(1,\ell)$. Then $\mu_\ell u$ and $\mu_\ell v$ belong to $\mathcal{P}(1,1)$. By \cite[Proposition~2.19]{reparam}, when $\mathcal{P}=\mathcal{M}$, there exist $\phi_u,\phi_v\in \mathcal{P}(1,1)$ such that $\mu_\ell u\phi_u = \mu_\ell v\phi_v$. When $\mathcal{P}=\mathcal{G}$, the same statement holds with $\phi_u=\id_1$ and $\phi_v=v^{-1} u$. In both cases, it implies that $u\phi_u = v\phi_v$, $\mu_\ell$ being invertible.
\epf

\bp \label{un}
Denote by $\mathcal{P}^1$ the full subcategory of $\mathcal{P}^{op}$ generated by $1$. Let $X$ be a $\mathcal{P}$-multipointed $d$-space. Let $(\alpha,\beta)\in X^0\p X^0$. The inclusion $j:\mathcal{P}^1\subset \mathcal{P}^{op}$ is final. Consequently, the topological space $\PA_{\alpha,\beta} X$ is also the quotient of $\PA^{top}_{\alpha,\beta} X$ by the reparametrization equivalence.
\ep

\bpf
Consider the comma category $(j\ddownarrow \ell)$ for a fixed $\ell>0$. It is nonempty because $\mathcal{P}(1,\ell)$ is nonempty. Let $u,v\in \mathcal{P}(1,\ell)$. Using Proposition~\ref{pre-un}, write $u\phi_u = v\phi_v$ for some $\phi_u,\phi_v\in \mathcal{P}(1,1)$. The equality $u\phi_u = v\phi_v$ means that there is the commutative diagram of $\mathcal{P}$
\[
\begin{tikzcd}[row sep=3em, column sep=3em]
	1 \arrow[d,->,"u"'] & 1 \arrow[l,"\phi_u"'] \arrow[r,,"\phi_v"] \arrow[d] & 1 \arrow[d,"v"']\\
	\ell \arrow[r,equal] & \ell \arrow[r,equal] & \ell
\end{tikzcd}
\]
In other terms, the comma category $(j\ddownarrow \ell)$ is connected. This implies that the inclusion functor $\mathcal{P}^1\subset \mathcal{P}^{op}$ is final in the sense of \cite[Section~IX.3]{MR1712872}. The proof is complete thanks to \cite[Theorem~IX.3.1]{MR1712872}.
\epf

\bd \label{path-with-length}
Let $X$ be a $\mathcal{P}$-multipointed $d$-space. The space $\PA^\ell_{\alpha,\beta}X$ is called the space of \textit{execution paths of length $\ell$} from $\alpha$ to $\beta$. Let \[\PA^\ell X = \coprod_{(\alpha,\beta) \in X^0\p X^0} \PA^\ell_{\alpha,\beta}X.\]
A map of multipointed $d$-spaces $f:X\to Y$ induces for each $\ell>0$ a continuous map $\PA^\ell f:\PA^\ell X \to \PA^\ell Y$ by composition by $f$ (in fact by $|f|$). The space $\PA^1X$ is also denoted by $\PA^{top}X$. 
\ed

\bp \label{variable-length}
Let $X$ be a $\mathcal{P}$-multipointed $d$-space. Let $\gamma_1$ and $\gamma_2$ be two execution paths of $X$ with $\gamma_1(1) = \gamma_2(0)$. Let $\ell_1,\ell_2>0$. Then $(\gamma_1\mu_{\ell_1} *  \gamma_2\mu_{\ell_1})\mu^{-1}_{\ell_1+\ell_2}$ is an execution path of $X$. 
\ep

\bpf
It is mutatis mutandis the proof of \cite[Proposition~4.10]{Moore2} after observing that $\mathcal{G}(1,1) \subset\mathcal{P}(1,1)$.
\epf

\bp \label{addlength} Let $X$ be a multipointed $d$-space. Let $\ell_1,\ell_2>0$. Let $(\alpha,\beta,\gamma)\in X^0\p X^0\p X^0$. 
The Moore composition of continuous maps yields a continuous maps \[\PA_{\alpha,\beta}^{\ell_1}X \p \PA_{\beta,\gamma}^{\ell_2}X \to \PA_{\alpha,\gamma}^{\ell_1+\ell_2}X. \]
\ep

\bpf It is a consequence of Proposition~\ref{variable-length}.
\epf

A possible reference for the notions of \textit{topological functor} and \textit{final structure} is \cite[Section~21]{topologicalcat}.

\bth \label{final-structure-revisited}
The mapping $\Omega:X\mapsto (|X|,X^0)$ induces a functor from $\ptop{\mathcal{P}}$ to the category $\mtop$ of multipointed spaces which is topological and fibre-small. The $\Omega$-final structure is the set of finite Moore compositions of the form $(f_1\gamma_1) * \dots * (f_n\gamma_n)$ such that $\gamma_i \in \PA^{\ell_i}X_i$ for all $1\leqslant i \leqslant n$ with $\sum_i \ell_i = 1$. 
\eth

Note that Theorem~\ref{final-structure-revisited} holds both by working with $\Delta$-generated spaces and with $\Delta$-Hausdorff $\Delta$-generated spaces. 

\bpf
The first statement is proved for $\mathcal{P}=\mathcal{G}$ in \cite[Proposition~6.5]{QHMmodel} using a description of the $\Omega$-initial structure which works for $\mathcal{P}=\mathcal{M}$ as well. The proof of the last statement is similar to the proof of \cite[Theorem~3.9]{Moore2}.
\epf

\bp \label{calcul-topology-glob} (\cite[Proposition~2.12]{Moore2} for $\mathcal{G}$ and $\mathcal{M}$)
Let $Z$ be a topological space of $\topspace$. Then there is the homeomorphism \[\PA_{0,1}^{top}\globP(Z)\iso \mathcal{P}(1,1)\p Z.\]
\ep

Because of the possible presence of stop intervals (see Definition~\ref{def_regular}) in the case $\mathcal{P}=\mathcal{M}$, the proof of Proposition~\ref{calcul-topology-glob} slightly differs from the proof of \cite[Proposition~2.12]{Moore2}. The latter, which uses an evaluation at $0.5$ to prove the continuity of $\pi$, is valid only for the case $\mathcal{P}=\mathcal{G}$. Indeed, in the case of $\mathcal{P}=\mathcal{M}$, for a given $\gamma\in \PA_{0,1}^{top}\globP(Z)$, $\gamma(0.5)$ can be equal to the initial or the final state of $\globP(Z)$.

\bpf
The set map $\Psi:\mathcal{P}(1,1)\p Z \to \PA_{0,1}^{top}\globP(Z)$ defined by $\Psi(\phi,z) = \delta_{z}\phi$ is continuous because the mapping $(t,\phi,z)\mapsto (z,\phi(t))$ from $[0,1]\p \mathcal{P}(1,1)\p Z$ to $|\globP(Z)|$ is continuous. It is a bijection since, by definition of $\globP(Z)$, the underlying set of $\PA_{0,1}^{top}\globP(Z)$ is equal to the underlying set of the space $\mathcal{P}(1,1)\p Z$. Consider the composite set map \[\pi:\PA_{0,1}^{top}\globP(Z) \longrightarrow \mathcal{P}(1,1)\p Z \longrightarrow Z\] which takes $\delta_{z}\phi$ to $\pr_2(\Psi^{-1}(\delta_{z}\phi))=z$ ($\pr_2$ is the projection on the second factor). Suppose that $\pi$ is not continuous. All involved topological spaces being sequential, there exist $z_\infty\in Z$, an open neighborhood $V$ of $z_\infty$ in $Z$, and a sequence $(\delta_{z_n}\phi_n)_{n\geq 0}$ which converges to the execution path $\delta_{z_\infty}\phi_\infty$ such that $z_n\in Z\backslash V$ for all $n\geq 0$. Choose $t_0\in [0,1]$ such that $\phi_\infty(t_0)\in ]0,1[$. The convergence for the compact-open topology implies the pointwise convergence. Thus the sequence $(\delta_{z_n}\phi_n(t_0))_{n\geq 0}$ of $|\globP(Z)|$ converges to $\delta_{z_\infty}\phi_\infty(t_0)$. This implies that there exists $N\geq 0$ such that for all $n\geq N$, $(z_n,\phi_n(t_0))\in |\globP(Z)|\backslash \{0,1\}$. By considering the image by the continuous projection (the left-hand term being equipped with the relative topology) \[|\globP(Z)|\backslash \{0,1\}\longrightarrow Z\] which is well-defined precisely because $0$ and $1$ are removed, we obtain that the sequence $(z_n)_{n\geq N}$ converges to $z_\infty$, and therefore that $z_\infty\in Z\backslash V$, the latter set being closed in $Z$: contradiction. This means that $\pi$ is continuous. The continuous map $Z\to \{0\}$ induces a continuous map 
\[
\begin{cases}
	& \PA_{0,1}^{top}\globP(Z)\longrightarrow \PA_{0,1}^{top}\globP(\{0\}) \iso \mathcal{P}(1,1) \\
	& \gamma \mapsto p.\gamma,
\end{cases}
\]
where $p:|\globP(Z)|\to [0,1]$ is the projection map. Therefore the set map 
\[
\begin{cases}
	& \Psi^{-1}:\PA_{0,1}^{top}\globP(Z)\longrightarrow \mathcal{P}(1,1) \p Z \\
	& \gamma \mapsto (p.\gamma,\pi(\gamma))
\end{cases}
\]
is continuous and $\Psi$ is a homeomorphism. 
\epf

\begin{cor} \label{quotient}
	Let $Z$ be a topological space of $\topspace$. Then there is the homeomorphism \[\PA_{0,1}\globP(Z)\iso Z.\]
\end{cor}

\bpf
There are the homeomorphisms 
\[
\liminj_{\mathcal{P}^1} \bigg(\mathcal{P}(1,1)\p Z\bigg) \iso \bigg(\liminj_{\mathcal{P}^1} \mathcal{P}(1,1)\bigg)\p Z \iso Z,
\]
the left-hand homeomorphism since $\topspace$ is cartesian closed and the right-hand homeomorphism by Proposition~\ref{pre-un}. The proof is complete thanks to Proposition~\ref{calcul-topology-glob} and Proposition~\ref{un}.
\epf

\bth \label{MG} 
Let $X$ be a $\mathcal{P}$-multipointed $d$-space. The following data assemble into a $\mathcal{P}$-flow denoted by $\moore^{\mathcal{P}}(X)$:
\begin{itemize}
	\item The set of states $X^0$ of $X$
	\item For all $\alpha,\beta\in X^0$, $\PA_{\alpha,\beta}\moore^{\mathcal{P}}(X) = \PA_{\alpha,\beta}^{\bullet}X$.
	\item For all $\alpha,\beta,\gamma\in X^0$ and all real numbers $\ell,\ell'>0$, the composition map is given by the map $*:\PA_{\alpha,\beta}^{\bullet} X \ot \PA_{\beta,\gamma}^{\bullet} X \to \PA_{\alpha,\gamma}^{\bullet} X$.
\end{itemize}
The mapping above induces a functor $\moore^{\mathcal{P}}:\ptop{\mathcal{P}}\to\dtopP$ which is a right adjoint.
\eth

\bpf
It is possible to construct a left adjoint by following step by step the method of \cite[Appendix~B]{Moore2}: the fact that all maps of $\mathcal{G}$ are invertible does not play any role at all.
\epf

\begin{nota}
	The left adjoint of $\moore^{\mathcal{P}}:\ptop{\mathcal{P}}\to\dtopP$ is denoted by \[\lmoore^{\mathcal{P}}:\dtopP \to \ptop{\mathcal{P}}.\] 
\end{nota}

\bp \label{precalcul} Let $X$ be a $\mathcal{P}$-multipointed $d$-space. Let $Z$ be a topological space. Then there is a bijection of sets \[\ptop{\mathcal{P}}(\globP(Z),X) \iso \coprod_{(\alpha,\beta)\in X^0\p X^0}\topspace(Z,\PA_{\alpha,\beta}^{top}X)\]
which is natural with respect to $Z$ and $X$.
\ep

\bpf
It is mutatis mutandis the proof of \cite[Proposition~4.13]{Moore2}.
\epf

\bp \label{calculM}
For all compact topological spaces $Z$, there are the natural isomorphisms $\moore^{\mathcal{P}}(\globP(Z)) \iso  \glob(\mathbb{F}_1^{\mathcal{P}^{op}}(Z))$ and $\mathbb{M}_!^{\mathcal{P}}(\glob(\mathbb{F}_1^{\mathcal{P}^{op}}(Z))) \iso  \globP(Z)$. 
\ep

\bpf
By definition of $\moore^{\mathcal{P}}$ and by Proposition~\ref{calcul-topology-glob}, the only nonempty path $\mathcal{P}$-space of $\moore^{\mathcal{P}}(\globP(Z))$ is $\PA_{0,1}\moore^{\mathcal{P}}(\globP(Z)) = \mathcal{P}(-,1)\p Z$: we obtain the first isomorphism. For any $\mathcal{P}$-multipointed $d$-space $X$, there is the sequence of natural bijections
\[\begin{aligned}
	\ptop{\mathcal{P}}\big(\mathbb{M}_!^{\mathcal{P}}(\glob(\mathbb{F}_1^{\mathcal{P}^{op}}(Z))),X \big)&\iso  \dtopP\big(\glob(\mathbb{F}_1^{\mathcal{P}^{op}}(Z)),\moore^{\mathcal{P}}X\big)\\
	&\iso \coprod_{(\alpha,\beta)\in X^0\p X^0} \topdgrP_0\big(\mathbb{F}_1^{\mathcal{P}^{op}}(Z),\PA_{\alpha,\beta}X\big)\\
	&\iso \coprod_{(\alpha,\beta)\in X^0\p X^0} \topspace(Z,\PA_{\alpha,\beta}^{1}X)\\
	&\iso \ptop{\mathcal{P}}(\globP(Z),X),
\end{aligned}\]
the first bijection by adjunction, the second bijection by \cite[Proposition~6.10]{Moore1}, the third bijection by Proposition~\ref{ev-adj} and the last bijection by Proposition~\ref{precalcul}. The proof of the second isomorphism is then complete thanks to the Yoneda lemma.
\epf

\section{Globular naturalization and carrier}
\label{core}

\begin{nota}
	Let $n\geq 1$. Denote by $\mathbf{D}^n = \{b\in \mathbb{R}^n, |b| \leqslant 1\}$ the $n$-dimensional disk, and by $\mathbf{S}^{n-1} = \{b\in \mathbb{R}^n, |b| = 1\}$ the $(n-1)$-dimensional sphere. By convention, let $\mathbf{D}^{0}=\{0\}$ and $\mathbf{S}^{-1}=\varnothing$.
\end{nota}

Every set $S$ can be viewed as a ``discrete'' $\mathcal{P}$-multipointed $d$-space $(S,S,\varnothing)$. The \textit{q-model structure} of $\mathcal{P}$-multipointed $d$-spaces is the unique combinatorial model structure such that 
\[\{\globP(\mathbf{S}^{n-1})\subset \globP(\mathbf{D}^{n}) \mid n\geq 0\} \cup \{C:\varnothing \to \{0\},R:\{0,1\} \to \{0\}\}\]
is the set of generating cofibrations, the maps between globes being induced by the closed inclusions $\mathbf{S}^{n-1}\subset \mathbf{D}^{n}$, and such that 
\[
\{\globP(\mathbf{D}^{n})\subset \globP(\mathbf{D}^{n+1}) \mid n\geq 0\}
\]
is the set of generating trivial cofibrations, the maps between globes being induced by the closed inclusions $(x_1,\dots,x_n)\mapsto (x_1,\dots,x_n,0)$. The weak equivalences are the maps of multipointed $d$-spaces $f:X\to Y$  inducing a bijection $f^0:X^0\iso Y^0$ and a weak homotopy equivalence $\PA^{{top}} f:\PA^{{top}} X \to \PA^{{top}} Y$ and the fibrations are the maps of multipointed $d$-spaces $f:X\to Y$  inducing a q-fibration $\PA^{{top}} f:\PA^{{top}} X \to \PA^{{top}} Y$ of topological spaces. A construction of this model structure is given in \cite[Theorem~6.16]{QHMmodel} for $\mathcal{P}=\mathcal{G}$. The argument works in the same way for $\mathcal{P}=\mathcal{M}$ because it relies on the use of the Quillen path object argument in \cite[Theorem~6.14]{QHMmodel} applied to the right adjoint from $\mathcal{P}$-multipointed $d$-spaces to topological graphs which forgets the composition and the reparametrization of execution paths. 
 
The space of execution paths of the cellular $\mathcal{P}$-multipointed $d$-spaces of the q-model structure are of particular interest. It is the purpose of this section to study them. Let $\lambda$ be an ordinal. We work with a colimit-preserving functor \[X:\lambda \longrightarrow \ptop{\mathcal{P}}\] such that
\begin{itemize}
	\item The $\mathcal{P}$-multipointed $d$-space $X_0$ is a set, in other terms $X_0=(X^0,X^0,\varnothing)$ for some set $X^0$.
	\item For all $\nu<\lambda$, there is a pushout diagram of $\mathcal{P}$-multipointed $d$-spaces 
	\[
	\begin{tikzcd}[row sep=3em, column sep=3em]
		\globP(\mathbf{S}^{n_\nu-1})\arrow[r,"g_\nu"] \arrow[d]  & X_\nu \arrow[d]\\
		\globP(\mathbf{D}^{n_\nu}) \arrow[r,"\widehat{g_\nu}"]  &  \cocartesian X_{\nu+1}
	\end{tikzcd}
	\]
	with $n_\nu \geq 0$. 
\end{itemize}

\begin{nota}
	Let $X_\lambda = \liminj_{\nu<\lambda} X_\nu$ be the transfinite composition.  
\end{nota}

\bd 
A \textit{cellular $\mathcal{P}$-multipointed $d$-space} is a $\mathcal{P}$-multipointed $d$-space of the form $X_\lambda$. 
\ed

For all $\nu\leqslant \lambda$, there is the equality $X_\nu^0=X^0$. Denote by \[c_\nu = |\globP(\mathbf{D}^{n_\nu})|\backslash |\globP(\mathbf{S}^{n_\nu-1})|\] the $\nu$-th cell of $X_\lambda$. It is called a \textit{globular cell}. Like in the usual setting of CW-complexes, $\widehat{g_\nu}$ induces a homeomorphism from $c_\nu$ to $\widehat{g_\nu}(c_\nu)$ equipped with the relative topology which will be therefore denoted in the same way. It also means that $\widehat{g_\nu}(c_\nu)$ equipped with the relative topology is $\Delta$-generated. The closure of $c_\nu$ in $|X_\lambda|$ is denoted by \[\widehat{c_{\nu}} = \widehat{g_{\nu}}(|\globP(\mathbf{D}^{n_\nu})|).\] 
It is a compact closed subset of $|X_\lambda|$. The boundary of $c_\nu$ in $|X_\lambda|$ is denoted by \[\de c_\nu = \widehat{g_{\nu}}(|\globP(\mathbf{S}^{n_\nu-1})|).\]
The state $\widehat{g_\nu}(0)\in X^0$ ($\widehat{g_\nu}(1)\in X^0$ resp.)  is called the \textit{initial (final resp.) state} of $c_\nu$. Note that they are not necessarily distinct. The integer $n_\nu+1$ is called the \textit{dimension} of the globular cell $c_\nu$. It is denoted by $\dim c_\nu$. The states of $X^0$ are also called the \textit{globular cells of dimension $0$}. 

As already noticed in \cite[Proposition~5.2]{Moore2}, the underlying space $|X_\lambda|$ is a cellular topological space. Proposition~\ref{normal} is required to be able to use \cite{reparam,reparam-fixed} for the Moore paths in $|X_\lambda|$.

\bp \label{normal}
The topological space $|X_\lambda|$ is Hausdorff. Let $K$ be a compact subspace of $|X_\lambda|$. Then $K$ intersects finitely many $c_\nu$.
\ep

\bpf The space $X_0$ is normal, being discrete. Adding one cell preserves normality by \cite[Proposition~1.1.2 (ii)]{MR1074175}. Assume that $\nu\leqslant \lambda$ is a limit ordinal and that each $X_\mu$ for $\mu<\nu$ is normal. We prove that $X_{\nu}$ is normal by an argument similar to the one of \cite[Proposition~A.5.1 (iv)]{MR1074175}. Each $X_\mu$ for $\mu\leqslant \lambda$ is $\Delta$-Hausdorff by \cite[Proposition~B.16]{leftproperflow}, and therefore has closed points. Hence $|X_\lambda|$ is Hausdorff. The last assertion is \cite[Proposition~5.5]{Moore2} whose proof mimicks \cite[Proposition~A.1]{MR1867354}.
\epf

\bp  \label{p1} \cite[Proposition~5.2]{Moore2}
The space $|X_\lambda|$ is a cellular space. It contains $X^0$ as a discrete closed subspace. For every $0\leqslant \nu_1 \leqslant \nu_2 \leqslant \lambda$, the continuous map $|X_{\nu_1}| \to |X_{\nu_2}|$ is a q-cofibration of spaces, and in particular a closed $T_1$-inclusion.
\ep

\bp \label{restriction_path} \cite[Proposition~5.3]{Moore2}
For all $0\leqslant \nu_1\leqslant \nu_2\leqslant\lambda$, there is the equality \[\PA^{top}X_{\nu_1} = \PA^{top} X_{\nu_2} \cap \ttop([0,1],|X_{\nu_1}|).\]
\ep

\bth  \label{topological-path-almost-accessible}
The composite functor 
\[
\lambda \stackrel{X} \longrightarrow \ptop{\mathcal{P}} \stackrel{\PA^{top}}\longrightarrow \topspace
\]
is colimit-preserving. In particular the continuous bijection \[\liminj (\PA^{top}.X) \longrightarrow \PA^{top} \liminj X\] is a homeomorphism. Moreover the topology of $\PA^{top} \liminj X$ is the final topology.
\eth

\bpf
It is mutatis mutandis the proof of \cite[Theorem~5.6]{Moore2} which relies on Proposition~\ref{restriction_path}.
\epf

\bth \label{cof-accessible} 
The composite functor 
\[
\lambda \stackrel{X} \longrightarrow \ptop{\mathcal{P}} \stackrel{\moore^{\mathcal{P}}}\longrightarrow \dtopP
\]
is colimit-preserving. In particular the natural map \[\liminj_{\nu<\lambda} \moore^{\mathcal{P}}(X_\nu) \longrightarrow \moore^{\mathcal{P}} X_\lambda\] is an isomorphism.
\eth

\bpf
The proof is mutatis mutandis the proof of \cite[Theorem~5.7]{Moore2} whose proof relies on Proposition~\ref{p1} and Theorem~\ref{topological-path-almost-accessible}.
\epf

\begin{nota}
	Let $c_\nu$ be a globular cell of $X_\lambda$. For all $z\in \mathbf{D}^{n_\nu}$ and all $\phi\in \mathcal{P}(1,1)$, the composite $\widehat{g_{\nu}}\delta_{z}\phi$ is an execution path of $X_\lambda$. When there is no ambiguity on the globular cell $c_\nu$ or no need to mention it, the execution path $\widehat{g_{\nu}}\delta_{z}\phi$ of $X_\lambda$ will be denoted by $\delta_{z}\phi$ (which is strictly speaking an abuse of notations already made in \cite{Moore2} to avoid overloading the notations).
\end{nota}

\bd \cite[Definition~1.1]{reparam} \label{def_regular} Let $\gamma:[0,L]\to |X_\lambda|$ be a Moore path. A \textit{stop interval} of $\gamma$ is an interval $[a,b]\subset [0,L]$ with $a<b$ such that the restriction $\gamma\rest_{[a,b]}$ is constant and such that $[a,b]$ is maximal for this property. The set of stop intervals of $\gamma$ is denoted by $\Delta_\gamma$. The path $\gamma$ is \textit{regular} if $\Delta_\gamma = \varnothing$ (no stop interval)~\footnote{Remember that by Definition~\ref{def_Moore_path}, all Moore paths of this paper are nonconstant.}. The Moore composition of two regular paths is regular. 
\ed 

The set $\Delta_\gamma$ contains only closed intervals since $|X_\lambda|$ is Hausdorff, Note that in the case $\mathcal{P}=\mathcal{G}$, all execution paths of $X_\lambda$ are regular by \cite[Proposition~5.13]{Moore2}. In the case  $\mathcal{P}=\mathcal{M}$, an execution path of the form $\widehat{g_{\nu}}\delta_{z}\phi$ with $\phi\in \mathcal{M}(1,1)\backslash \mathcal{G}(1,1)$ is not regular. Proposition~\ref{example_reg} and Proposition~\ref{local-inj-reg} play a key role in the sequel.

\bp \label{example_reg}
Consider a globular cell $c_\nu$ of $X_\lambda$. Let $z\in \mathbf{D}^{n_\nu}\backslash \mathbf{S}^{n_\nu-1}$. The execution path $\widehat{g_{\nu}}\delta_{z}$ is regular.
\ep

Note that $\widehat{g_{\nu}}\delta_{z}$ is not necessarily regular if $z\in \mathbf{S}^{n_\nu-1}$.

\bpf
It is a consequence of the facts that $\widehat{g_{\nu}}$ induces a homeomorphism from $c_\nu$ to $\widehat{g_{\nu}}(c_\nu)$ and that $\delta_z$ is a regular execution path of $\globP(\mathbf{D}^{n_\nu})$. 
\epf

\bp \label{local-inj-reg}
Let $\gamma:[0,L]\to |X_\lambda|$ be a Moore path. Suppose that $\gamma$ is regular and that there exist $\eta_1,\eta_2\in \mathcal{P}(\ell,L)$ such that $\gamma\eta_1=\gamma\eta_2$. Then $\eta_1=\eta_2$. 
\ep

\bpf Note that $\eta_1,\eta_2\in \mathcal{P}(\ell,L) \subset \mathcal{M}(\ell,L)$. From $\gamma\eta_1=\gamma\eta_2$, we obtain $\gamma\mu_L^{-1}\mu_L\eta_1\mu_\ell^{-1}=\gamma\mu_L^{-1}\mu_L\eta_2\mu_\ell^{-1}$. The Moore path $\gamma\mu_L^{-1}:[0,1]\to |X_\lambda|$ is regular, $\mu_L^{-1}$ being a homeomorphism. Using \cite[Lemma~3.9]{reparam}, we deduce that $\mu_L\eta_1\mu_\ell^{-1} = \mu_L\eta_2\mu_\ell^{-1}$ and therefore that $\eta_1=\eta_2$.
\epf

\bth \label{normal-form}
Let $\gamma$ be an execution path of $X_\lambda$. It can be decomposed as a Moore composition 
\[
\gamma = (\widehat{g_{\nu_1}}\delta_{z_1}\phi_1\mu_{\ell_1})*\dots * (\widehat{g_{\nu_n}}\delta_{z_n}\phi_n\mu_{\ell_n})
\]
with $n\geq 1$, $\nu_i<\lambda$ and $z_i\in \mathbf{D}^{n_{\nu_i}}\backslash \mathbf{S}^{n_{\nu_i}-1}$ and $\phi_i\in \mathcal{P}(1,1)$ for all $i\in \{1,\dots,n\}$ and $\ell_1 + \dots + \ell_n= 1$. Consider a second decomposition 
\[
\gamma  = (\widehat{g_{\nu'_1}}\delta_{z'_1}\phi'_1\mu_{\ell'_1})*\dots * (\widehat{g_{\nu'_{n'}}}\delta_{z'_{n'}}\phi'_n\mu_{\ell'_{n'}}).
\]
Then $n=n'$, $\nu_i=\nu'_i$ and $z_i=z'_i$ for all $1\leqslant i \leqslant n$ and 
\[
(\phi_1\mu_{\ell_1})\ot\dots \ot (\phi_n\mu_{\ell_n}) = (\phi'_1\mu_{\ell'_1})\ot \dots \ot (\phi'_n\mu_{\ell'_{n}})\in \mathcal{P}(1,n).
\]
\eth

\bpf
By Theorem~\ref{final-structure-revisited}, every execution path $\gamma$ from $\alpha$ to $\beta$ of $X_\lambda$ is of the form a Moore composition $\gamma = (\widehat{g_{\nu_1}}\delta_{z_1}\phi_1\mu_{\ell_1})*\dots * (\widehat{g_{\nu_n}}\delta_{z_n}\phi_n\mu_{\ell_n})$ with $n\geq 1$, $\nu_i<\lambda$ and $z_i\in \mathbf{D}^{n_{\nu_i}}\backslash \mathbf{S}^{n_{\nu_i}-1}$ and $\phi_i\in \mathcal{P}(1,1)$ for all $i\in \{1,\dots,n\}$ and $\ell_1 + \dots + \ell_n= 1$. Consider a second decomposition $\gamma  = (\widehat{g_{\nu'_1}}\delta_{z'_1}\phi'_1\mu_{\ell'_1})*\dots * (\widehat{g_{\nu'_{n'}}}\delta_{z'_{n'}}\phi'_n\mu_{\ell'_{n'}})$. Then $\gamma$ is the Moore composition of a Moore path going from $\widehat{g_{\nu_1}}(0)$ to $\widehat{g_{\nu_1}}(1)$ in the globular cell $c_{\nu_1}$ followed by a Moore path going from $\widehat{g_{\nu_2}}(0)$ to $\widehat{g_{\nu_2}}(1)$ in the globular cell $c_{\nu_2}$ etc... until a Moore path going from $\widehat{g_{\nu_n}}(0)$ to $\widehat{g_{\nu_n}}(1)$ in the globular cell $c_{\nu_n}$. And $\gamma$ is also the Moore composition of a Moore path going from $\widehat{g_{\nu'_1}}(0)$ to $\widehat{g_{\nu'_1}}(1)$ in the globular cell $c_{\nu'_1}$ followed by a Moore path going from $\widehat{g_{\nu'_2}}(0)$ to $\widehat{g_{\nu'_2}}(1)$ in the globular cell $c_{\nu'_2}$ etc... until a Moore path going from $\widehat{g_{\nu'_{n'}}}(0)$ to $\widehat{g_{\nu'_{n'}}}(1)$ in the globular cell $c_{\nu'_{n'}}$. From the set bijection 
\[
|X_\lambda| = X^0 \sqcup \coprod_{\nu<\lambda} c_\nu,
\]
we deduce that $n=n'$, $\nu_i=\nu'_i$ and $z_i=z'_i$ for $1\leqslant i \leqslant n$. By \cite[Proposition~3.4]{Moore2}, we also have 
\begin{multline*}
	\bigg((\widehat{g_{\nu_1}}\delta_{z_1})*\dots * (\widehat{g_{\nu_n}}\delta_{z_n})\bigg)\bigg((\phi_1\mu_{\ell_1})\ot\dots \ot (\phi_n\mu_{\ell_n})\bigg) \\ = \bigg((\widehat{g_{\nu'_1}}\delta_{z'_1})*\dots * (\widehat{g_{\nu'_{n'}}}\delta_{z'_{n'}})\bigg) \bigg((\phi'_1\mu_{\ell'_1})\ot \dots \ot (\phi'_n\mu_{\ell'_{n}})\bigg).
\end{multline*}
Since $(\widehat{g_{\nu_1}}\delta_{z_1})*\dots * (\widehat{g_{\nu_n}}\delta_{z_n})=(\widehat{g_{\nu'_1}}\delta_{z'_1})*\dots * (\widehat{g_{\nu'_{n'}}}\delta_{z'_{n'}})$ is a regular Moore path, being a Moore composition of regular Moore paths by Proposition~\ref{example_reg}, we deduce by Proposition~\ref{local-inj-reg} the equality \[(\phi_1\mu_{\ell_1})\ot\dots \ot (\phi_n\mu_{\ell_n}) = (\phi'_1\mu_{\ell'_1})\ot \dots \ot (\phi'_n\mu_{\ell'_{n}}).\]
\epf

Theorem~\ref{normal-form} is a modification of \cite[Theorem~5.9]{Moore2} which is valid both for $\mathcal{G}$ and $\mathcal{M}$. Let $\Psi=(\phi_1\mu_{\ell_1})\ot\dots \ot (\phi_n\mu_{\ell_n})= (\phi'_1\mu_{\ell'_1})\ot \dots \ot (\phi'_n\mu_{\ell'_{n}})$. If $\mathcal{P}=\mathcal{G}$, then $\ell_1+\dots + \ell_i = \Psi^{-1}(i)$ for $1\leqslant i \leqslant n$. This implies that $\ell_i=\ell'_i$ for $1\leqslant i \leqslant n$, and by Proposition~\ref{decomposition-tenseur}, we deduce that $\phi_i=\phi'_i$ for $1\leqslant i \leqslant n$, which implies \cite[Theorem~5.9]{Moore2}.

\bd \label{def-regularization}
With the notations of Theorem~\ref{normal-form}, the regular Moore path
\[
\natgl(\gamma) = (\widehat{g_{\nu_1}}\delta_{z_1})*\dots * (\widehat{g_{\nu_n}}\delta_{z_n})
\]
is called the \textit{(globular) naturalization} of $\gamma$. The sequence of globular cells \[\carrier(\gamma)=[c_{\nu_1},\dots,c_{\nu_n}]\] is called the \textit{carrier} of $\gamma$. The integer $n$ is called the \textit{length} of the carrier. It is also called the \textit{natural length} of $\gamma$. 
\ed

The globular naturalization of the Moore composition of two execution paths is the Moore composition of the globular naturalizations. 

\bp \label{naturalization-is-execution-path}
Let $\gamma$ be an execution path of $X_\lambda$ of natural length $n$. Then the regular Moore path $\natgl(\gamma)$ is an execution path of $X_\lambda$ of length $n$.
\ep

\bpf
It is a consequence of Proposition~\ref{addlength}.
\epf

\bd \label{minimal} An execution path $\gamma$ of $X_\lambda$ is minimal~\footnote{It is not exactly the definition chosen in \cite{Moore2}. This one makes sense only for cellular $\mathcal{P}$-multipointed $d$-spaces.} if  $\gamma=\widehat{g_{\nu}}\delta_{z}\phi$ for some $\nu<\lambda$, some $z\in \mathbf{D}^{n_\nu} \backslash \mathbf{S}^{n_\nu-1}$ and some $\phi\in \mathcal{P}(1,1)$. 
\ed

\bth \label{calcul_final_structure} Let $0\leqslant \nu<\lambda$. 
Then every execution path of $X_{\nu+1}$ can be written as a finite Moore composition 
\[(f_1\gamma_1\mu_{\ell_1}) * \dots * (f_n\gamma_n\mu_{\ell_n})\] with $n\geq 1$ such that 
\begin{enumerate}[leftmargin=*]
	\item $\sum_i \ell_i = 1$.
	\item $f_i = g_\nu$ and $\gamma_i$ is an execution path of $X_\nu$ or $f_i=\widehat{g_\nu}$ and $\gamma_i=\delta_{z_i}\phi_i$ with $z_i\in \mathbf{D}^{n_\nu}\backslash\mathbf{S}^{n_\nu-1}$ and some $\phi_i\in \mathcal{P}(1,1)$.
	\item for all $1\leqslant i <n$, either $f_i\gamma_i$ or $f_{i+1}\gamma_{i+1}$ (or both) is (are) of the form $\widehat{g_\nu}\delta_{z}\phi$ for some $z\in \mathbf{D}^{n_\nu}\backslash\mathbf{S}^{n_\nu-1}$ and some $\phi\in \mathcal{P}(1,1)$: intuitively, there is no possible simplification using the Moore composition inside $X_\nu$. 
\end{enumerate}
If there is another finite Moore composition 
\[(f'_1\gamma'_1\mu_{\ell'_1}) * \dots * (f'_n\gamma'_n\mu_{\ell'_{n'}})\] with $n'\geq 1$ satisfying the same properties, then $n=n'$, for all $i\in \{1,\dots,n\}$ one has $f_i=f'_i$, $r_i=\natgl(\gamma_i)=\natgl(\gamma'_i)$ and finally 
\[
(\phi_1\mu_{\ell_1})\ot \dots \ot (\phi_n\mu_{\ell_n}) = (\phi'_1\mu_{\ell'_1})\ot \dots \ot (\phi'_n\mu_{\ell'_n})
\]
with $\gamma_i=r_i\phi_i$ and $\gamma'_i=r_i\phi'_i$ for all $i\in \{1,\dots,n\}$. 
\eth

\bpf The existence of the finite Moore composition is a consequence of Theorem~\ref{final-structure-revisited}. Let $i\in \{1,\dots,n\}$. If $f_i=f$, then $f_i$ is one-to-one by Proposition~\ref{p1}. Thus, the execution path $f\natgl(\gamma_i)$ is regular. Besides, $f\natgl(\gamma_i)$ is in this case the regular naturalization of $f\gamma_i$. If $f_i=\widehat{g_\nu}$ and $\gamma_i=\delta_{z_i}\phi_i$, then $f_i\delta_{z_i}$ is regular by Proposition~\ref{example_reg}. And by definition of the globular naturalization, $f_i\delta_{z_i} = \natgl(f_i\delta_{z_i})$. Therefore we obtain 
\[
(f_1\natgl(\gamma_1)) * \dots * (f_n\natgl(\gamma_n)) = (f'_1\natgl(\gamma'_1)) * \dots * (f'_n\natgl(\gamma'_n)) \in \PA^n X_\lambda.
\]
By definition of the Moore composition, it implies that $\natgl(\gamma_i) = \natgl(\gamma'_i)$ for all $i\in \{1,\dots,n\}$. The last equality is a consequence of Proposition~\ref{local-inj-reg}.
\epf

Theorem~\ref{calcul_final_structure} is a replacement for \cite[Theorem~5.20]{Moore2} which is valid both for $\mathcal{G}$ and $\mathcal{M}$. Let $\Psi=(\phi_1\mu_{\ell_1})\ot\dots \ot (\phi_n\mu_{\ell_n})= (\phi'_1\mu_{\ell'_1})\ot \dots \ot (\phi'_n\mu_{\ell'_{n}})$. If $\mathcal{P}=\mathcal{G}$, then $\ell_1+\dots + \ell_i = \Psi^{-1}(i)$ for $1\leqslant i \leqslant n$. This implies that $\ell_i=\ell'_i$ for $1\leqslant i \leqslant n$, and by Proposition~\ref{decomposition-tenseur}, we deduce that $\phi_i=\phi'_i$ for $1\leqslant i \leqslant n$, which implies \cite[Theorem~5.20]{Moore2}.

\bd Let $c_\nu$ be a globular cell of the cellular $\mathcal{P}$-multipointed $d$-space $X_\lambda$ with $\nu<\lambda$ and $\dim(c_\nu)\geq 1$. Let $0<h<1$. Let 
\[
\widehat{c_{\nu}}[h] = \bigg\{\widehat{g_{\nu}}(z,h)\mid (z,h)\in |\globP(\mathbf{D}^{n_{\nu}})|\bigg\}.
\]
It is called an \textit{achronal slice} of the globular cell $c_\nu$. 
\ed

\cite[Proposition~5.17]{Moore2} claims that, for any globular cell $c_\nu$ of any cellular $\mathcal{G}$-multipointed $d$-space $X_\lambda$ with $\dim(c_\nu)\geq 1$, there exists $b\in ]0,1[$ such that for all $h\in ]0,b]$, one has  $\widehat{c_{\nu}}[h]\cap X^0=\varnothing$. This implies that there exists $h\in ]0,1[$ such that $\widehat{c_{\nu}}[h]\cap X^0=\varnothing$. In plain English, this means that there is an achronal slice of the globular cell $c_\nu$ which does not intersect $X^0$. It is the key fact to prove \cite[Theorem~5.18]{Moore2}, and then to deduce \cite[Theorem~5.19]{Moore2}. Proposition~\ref{impossible} proves that \cite[Proposition~5.17]{Moore2} is false for $\mathcal{P}=\mathcal{M}$.

\bp \label{impossible}
There exists a cellular $\mathcal{M}$-multipointed $d$-space $X_\lambda$ and a globular cell $c_\nu$ with $\nu<\lambda$ and $\dim(c_\nu)\geq 1$ such that for all $h\in ]0,1[$, $\widehat{c_{\nu}}[h]\cap X^0\neq \varnothing$.
\ep

\bpf
Consider the continuous map $\Psi:[-1,1]\p [0,1]\to [0,1]$ defined by 
\[
\Psi:(x,t)\mapsto 
\begin{cases}
	0 & \hbox{ if }  0\leqslant t \leqslant \frac{2+x}{4}\\
	4t-2-x & \hbox{ if } \frac{2+x}{4}\leqslant t \leqslant \frac{3+x}{4}\\
	1 & \hbox{ if } \frac{3+x}{4} \leqslant t \leqslant 1
\end{cases}
\]
The continuous map \[f:((x,y),t) \mapsto ((x,y),\Psi(x,t))\] from $\mathbf{S}^1\p [0,1]$ to $|\globM(\mathbf{S}^1)|$ induces a continuous map from $|\globM(\mathbf{S}^1)|$ to itself since $f((x,y),0)=((x,y),\Psi(x,0)) = 0$ and $f((x,y),1)=((x,y),\Psi(x,1)) = 1$ for all $(x,y) \in \mathbf{S}^1$ by definition of $\Psi$. An execution path of the form $\delta_{(x,y)}\phi$ with $\phi\in \mathcal{M}(1,1)$ is taken by this continuous map to the continuous path $\delta_{(x,y)}\Psi(x,\phi(t))$. The latter is an execution path of $|\globM(\mathbf{S}^1)|$ since $\Psi(x,-)\in \mathcal{M}(1,1)$ for all $x\in [-1,1]$. Consequently, we obtain a map of $\mathcal{M}$-multipointed $d$-spaces $f:\globM(\mathbf{S}^1)\to\globM(\mathbf{S}^1)$. Consider the $\mathcal{M}$-multipointed $d$-space $X$ defined by the pushout diagram of $\mathcal{M}$-multipointed $d$-spaces
\[
\begin{tikzcd}[row sep=3em, column sep=3em]
	\globM(\mathbf{S}^1) \arrow[d] \arrow[r,"f"] & \globM(\mathbf{S}^1) \arrow[d] \\
	\globM(\mathbf{D}^2) \arrow[r,"\widehat{f}"] & \cocartesian X
\end{tikzcd}
\]
where the left vertical map is induced by the inclusion $\mathbf{S}^1 \subset \mathbf{D}^2$. The $\mathcal{M}$-multipointed space $X$ is cellular, $\globM(\mathbf{S}^1)$ being cellular. Let $h\in ]0,1[$. Consider the achronal slice \[\widehat{f}[h] = \{\widehat{f}(z,h)\mid (z,h)\in |\globM(\mathbf{D}^2)|\}.\] One has \[\widehat{f}((0,1),h)= f((0,1),h) = ((0,1),\Psi(0,h)),\]
the first equality because the square above is commutative and the second equality by definition of $f$. This implies that $\widehat{f}((0,1),h)=0$ when $h\leqslant 1/2$ by definition of $\Psi$. Similarly, there are the equalities \[\widehat{f}((-1,0),h)= f((-1,0),h) = ((-1,0),\Psi(-1,h)).\] This implies that $\widehat{f}((-1,0),h)=1$ when $h\geq 1/2$ by definition of $\Psi$. 

We deduce that for all $h\in ]0,1[$, $\widehat{f}[h]\cap X^0\neq \varnothing$.
\epf

In this paper, we prove Theorem~\ref{diagonal-execution} instead. It enables us to deduce both \cite[Theorem~5.18]{Moore2} and \cite[Theorem~5.19]{Moore2} in a different way for the two reparametrization categories $\mathcal{G}$ and $\mathcal{M}$.

\begin{nota}
	Let $(\alpha,\beta)\in X^0\p X^0$. Denote by $(\PA_{\alpha,\beta}^{top}X_\lambda)_{co}$ the set $\PA_{\alpha,\beta}^{top}X_\lambda$ equipped with the compact-open topology.
\end{nota}

\bth \label{diagonal-execution} (replacement for \cite[Proposition~5.17]{Moore2})
Let $(\alpha,\beta)\in X^0\p X^0$. Let $(\gamma_n)_{n\geq 0}$ be a sequence of $(\PA_{\alpha,\beta}^{top}X_\lambda)_{co}$ which converges to $\gamma_\infty$. Then the set $\{\carrier(\gamma_n)\mid n\geq 0\}$ is finite.
\eth

\bpf
Consider the one-point compactification $\overline{\mathbb{N}}=\mathbb{N}\cup \{\infty\}$ of the discrete space of integers $\mathbb{N}$. Note that $\overline{\mathbb{N}}$ is not $\Delta$-generated, its $\Delta$-kelleyfication being discrete. The converging sequence $(\gamma_n)_{n\geq 0}$ gives rise to a continuous map \[\psi:\overline{\mathbb{N}} \longrightarrow (\PA_{\alpha,\beta}^{top}X_\lambda)_{co} \subset \ttop_{co}([0,1],|X_\lambda|)\] where $\ttop_{co}([0,1],|X_\lambda|)$ is the set of continuous maps from $[0,1]$ to $|X_\lambda|$ equipped with the compact-open topology. Since $[0,1]$ is locally compact, it is exponential in the category of \textit{general} topological spaces by \cite[Proposition~7.1.5]{Borceux2}. We obtain a continuous map \[\widehat{\psi}:\overline{\mathbb{N}}\p_{gen} [0,1] \longrightarrow |X_\lambda|\] where $\p_{gen}$ is the binary product in the category of general topological spaces. Since $\overline{\mathbb{N}}\p_{gen} [0,1]$ is compact by Tychonoff, the subset $\widehat{\psi}(\overline{\mathbb{N}}\p_{gen} [0,1])$ is compact and closed in $|X_\lambda|$, the latter being Hausdorff by Proposition~\ref{normal}. The subset $\widehat{\psi}(\overline{\mathbb{N}}\p_{gen} [0,1])$ therefore intersects finitely many globular cells $\{c_{\nu_j}\mid j\in J\}$ by Proposition~\ref{normal}. Suppose that the set $\{\carrier(\gamma_n)\mid n\geq 0\}$ is infinite. This implies that the sequence of lengths of $\carrier(\gamma_n)$ for $n\geq 0$ is not bounded, the set $J$ being finite. By extracting a subsequence, one can suppose that the sequence of lengths is strictly increasing. Each infinite sequence of $\{c_{\nu_j}\mid j\in J\}$ has a constant infinite subsequence since $J$ is finite. Therefore by a Cantor diagonalization argument, one can suppose that there exists a sequence $(c_{\nu_{j_n}})_{n\geq 0}$ of $\{c_{\nu_j}\mid j\in J\}$ such that for all $n\geq 0$, there is the equality \[\carrier(\gamma_n) = [c_{\nu_{j_0}},\dots,c_{\nu_{j_{i_n}}}]\] for some strictly increasing sequence of integers $(i_n)_{n\geq 0}$. This means that the execution path $\gamma_n$ is the composition of an execution path whose image is included in $\{\widehat{g_{\nu_{j_0}}}(0),\widehat{g_{\nu_{j_0}}}(1)\} \cup c_{\nu_{j_0}}$ from $\widehat{g_{\nu_{j_0}}}(0)$ to $\widehat{g_{\nu_{j_0}}}(1)$ followed by an execution path whose image is included in  $\{\widehat{g_{\nu_{j_1}}}(0)\cup \widehat{g_{\nu_{j_1}}}(1)\} \cup c_{\nu_{j_1}}$ from $\widehat{g_{\nu_{j_0}}}(1)=\widehat{g_{\nu_{j_1}}}(0)$ to $\widehat{g_{\nu_{j_1}}}(1)$ etc... until an execution path whose image is included in $\{\widehat{g_{\nu_{j_{i_n}}}}(0)\cup \widehat{g_{\nu_{j_{i_n}}}}(1)\} \cup c_{\nu_{j_{i_n}}}$. The sequence of execution paths $(\gamma_n)_{n\geq 0}$ converges pointwise to $\gamma_\infty$ because it converges for the compact-open topology by hypothesis. The subset $\widehat{c_{\nu}}$ being closed in $|X_\lambda|$ for all $\nu<\lambda$, the sequence of globular cells $\carrier(\gamma_\infty)$ consists of a concatenation of sequences $S_{\nu_{j_n}}$ for $n\geq 0$ where either $S_{\nu_{j_n}} = [c_{\nu_{j_n}}]$ or $S_{\nu_{j_n}}$ is a nonempty finite sequence of globular cells intersecting $\de c_{\nu_{j_n}}$. This implies that the sequence of globular cells $\carrier(\gamma_\infty)$ is infinite: contradiction.
\epf

\begin{cor} \label{bounded0} (\cite[Theorem~5.18]{Moore2} for $\mathcal{G}$ and $\mathcal{M}$)
Let $\gamma_\infty$ be an execution path of $X_\lambda$. Let $\nu_0<\lambda$. There exists an open neighborhood $\Omega$ of $\gamma_\infty$ in $\PA^{top}X_\lambda$ such that for all execution paths $\gamma\in \Omega$, the number of copies of $c_{\nu_0}$ in the carrier of $\gamma$ does not exceed the length of the carrier of $\gamma_\infty$.
\end{cor}

\bpf
Let $\Omega_{\nu_0}$ be the set of execution paths $\gamma$ such that the number of copies of $c_{\nu_0}$ in the carrier of $\gamma$ does not exceed the length of the carrier of $\gamma_\infty$. Suppose that $\gamma_\infty$ is not in the interior of $\Omega_{\nu_0}$. Since $\PA^{top}X_\lambda$ is sequential, being $\Delta$-generated, there exists a sequence $(\gamma_n)_{n\geq 0}$ of the complement of $\Omega_{\nu_0}$ converging to $\gamma_\infty$. By Theorem~\ref{diagonal-execution}, the set $\{\carrier(\gamma_n)\mid n\geq 0\}$ is finite. Thus by extracting a subsequence, one can suppose that the sequence of carriers $(\carrier(\gamma_n))_{n\geq 0}$ is constant, write \[\carrier(\gamma_n) = [c_{\nu_{1}},\dots,c_{\nu_{N}}]\] for all $n\geq 0$. The integer $N$ is strictly greater than the length of $\carrier(\gamma_\infty)$ since $[c_{\nu_{1}},\dots,c_{\nu_{N}}]$ contains strictly more copies of $c_{\nu_0}$ than the length of $\carrier(\gamma_\infty)$ by definition of $\Omega_{\nu_0}$. The sequence $(\gamma_n)_{n\geq 0}$ converges also pointwise to $\gamma_\infty$. Thus, $\carrier(\gamma_\infty)$ consists of a concatenation of sequences $S_{\nu_{n}}$ for $1\leqslant n\leqslant N$ where either $S_{\nu_{n}} = [c_{\nu_{n}}]$ or $S_{\nu_{n}}$ is a nonempty finite sequence of globular cells intersecting $\de c_{\nu_{n}}$. This implies that the length of $\carrier(\gamma_\infty)$ is strictly greater than itself: contradiction. We deduce that $\gamma_\infty$ is in the interior of $\Omega_{\nu_0}$. Hence the existence of the open neighborhood.
\epf

Corollary~\ref{bounded0} proves the existence of an open neighborhood $\Omega$ in the $\Delta$-kelleyfication of the compact-open topology. The latter topology contains more open subsets than the compact-open topology. The proof of \cite[Theorem~5.18]{Moore2} implies that $\Omega$ can even be taken in the compact-open topology when $\mathcal{P}=\mathcal{G}$ by using the formula of \cite[page~197]{Moore2}. Because of Proposition~\ref{impossible}, the formula of \cite[page~197]{Moore2} does not necessarily provide anymore an open of the compact-open topology for the case $\mathcal{P} = \mathcal{M}$. We do not know whether the formula of \cite[page~197]{Moore2} could provide in some cases an open subset of the $\Delta$-kelleyfication of the compact-open topology, which could supply an alternative proof of Corollary~\ref{bounded0}.

\begin{cor} \label{bounded} (\cite[Theorem~5.19]{Moore2} for $\mathcal{G}$ and $\mathcal{M}$)
Let $(\gamma_k)_{k\geq 0}$ be a sequence of execution paths of $X_\lambda$ which converges in $\PA^{top}X_\lambda$. Let $c_{\nu_0}$ be a globular cell of $X_\lambda$. Let $i_k$ be the number of times that $c_{\nu_0}$ appears in $\carrier(\gamma_k)$. Then the sequence of integers $(i_k)_{k\geq 0}$ is bounded.
\end{cor}

\bpf
The sequence $(\gamma_k)_{k\geq 0}$ converges in $(\PA^{top}X_\lambda)_{co}$ as well because of the continuous map $\PA^{top}X_\lambda\to (\PA^{top}X_\lambda)_{co}$. Thus the set $\{\carrier(\gamma_n)\mid n\geq 0\}$ is finite by Theorem~\ref{diagonal-execution} and, therefore, the sequence of integers $(i_k)_{k\geq 0}$ is bounded.
\epf

Theorem~\ref{carrier-finite-on-compact} is a much better statement that will be the replacement in this paper of \cite[Theorem~5.19]{Moore2}. 

\bth \label{carrier-finite-on-compact}
Let $(\alpha,\beta)\in X^0\p X^0$. Let $\psi:[0,1] \to \PA_{\alpha,\beta}^{top}X_\lambda$ be a continuous map. Then the set $\{\carrier(\psi(u)) \mid u\in [0,1]\}$ is finite.
\eth

\bpf
Suppose that the set $\{\carrier(\psi(u)) \mid u\in [0,1]\}$ is infinite. Then there exists a sequence $(t_n)_{n\geq 0}$ of $[0,1]$ such that 
\[
\forall m,n\geq 0, m\neq n \Rightarrow \carrier(\psi(t_m))\neq \carrier(\psi(t_n)). 
\]
In particular, this means that $t_m=t_n$ implies $m=n$ for all $m,n\geq 0$. By extracting a subsequence, one can suppose that the sequence $(t_n)_{n\geq 0}$ converges to some $t_\infty\in [0,1]$. And the above condition ensures that the set of carriers $\{\carrier(\psi(t_n)) \mid n\geq 0\}$ is still infinite. Since $\psi$ is continuous, the sequence of execution paths $(\psi(t_n))_{n\geq 0}$ converges to $\psi(t_\infty)$ in $\PA_{\alpha,\beta}^{top}X_\lambda$, and therefore in $(\PA_{\alpha,\beta}^{top}X_\lambda)_{co}$ because of the continuous map $\PA_{\alpha,\beta}^{top}X_\lambda\to(\PA_{\alpha,\beta}^{top}X_\lambda)_{co}$. That contradicts Theorem~\ref{diagonal-execution}.
\epf

Theorem~\ref{carrier-finite-on-compact} enables us to understand the difference between Raussen's naturalization of \cite[Definition~2.14]{MR2521708} and the globular naturalization of Definition~\ref{def-regularization}.

\begin{cor} \label{diff-glob-cube}
	Let $(\alpha,\beta)\in X^0\p X^0$. Let $\psi:[0,1] \to \PA_{\alpha,\beta}^{top}X_\lambda$ be a continuous map. Then the set of natural lengths of $\psi(u)$ for $u$ running over $[0,1]$ is bounded.
\end{cor}

\bpf
It is due to the fact that the set $\{\carrier(\psi(u)) \mid u\in [0,1]\}$ is finite by Theorem~\ref{carrier-finite-on-compact}.
\epf

The natural lengths of $\psi(u)$ for $u$ running over $[0,1]$ have no reason to be constant. Consider a pushout of $\mathcal{P}$-multipointed $d$-spaces 
\[
\begin{tikzcd}[row sep=3em, column sep=3em]
	\globP(\mathbf{S}^0) \arrow[r] \arrow[d] & \vI^\mathcal{P} * \vI^\mathcal{P} \arrow[d] \\
	\globP(\mathbf{D}^1) \arrow[r] & \cocartesian X
\end{tikzcd}
\]
where $\vI^\mathcal{P} * \vI^\mathcal{P}$ means that the final state of the left copy of the directed segment is identified with the initial state of the right copy of the directed segment (see Notation~\ref{globchain}) and where the top horizontal map (it is not unique) takes the initial (final resp.) state of $\globP(\mathbf{S}^0)$ to the initial (final resp.) state of $\vI^\mathcal{P} * \vI^\mathcal{P}$. Then the natural length of an execution path of $X$ going from the initial to the final state is $2$ on the boundary of $\globP(\mathbf{D}^1)$ and $1$ inside $\globP(\mathbf{D}^1)$.

Corollary~\ref{diff-glob-cube} shows the difference of behavior between the natural length in the globular setting and the length of the naturalization of a directed path between two vertices in the geometric realization of a precubical set. Indeed, the latter is constant by continuous deformation preserving the extremities \cite[Section~2.2.1]{MR2521708} \cite[Proposition~2.2]{Raussen2012}. See also \cite[Proposition~4.5]{RegularMoore}.

\section{Chain of globes}
\label{core2}

\begin{nota} \label{globchain}
	Let $Z_1,\dots,Z_p$ be $p$ nonempty topological spaces with $p\geq 1$. Consider the $\mathcal{P}$-multipointed $d$-space \[X=\globP(Z_1)* \dots *\globP(Z_p)\] with $p\geq 1$ where the $*$ means that the final state of a globe is identified with the initial state of the next one by reading from the left to the right. Let $\{\alpha_0,\alpha_1,\dots,\alpha_p\}$ be the set of states such that the canonical map $\globP(Z_i)\to X$ takes the initial state $0$ of $\globP(Z_i)$ to $\alpha_{i-1}$ and the final state $1$ of $\globP(Z_i)$ to $\alpha_{i}$.
\end{nota}

\begin{nota} \label{Xc}
	Each carrier $\underline{c}=[c_{\nu_1},\dots,c_{\nu_n}]$ gives rise to a map of $\mathcal{P}$-multipointed $d$-spaces from a chain of globes to $X_\lambda$
\[
\widehat{g_{\underline{c}}}:\globP(\mathbf{D}^{n_{\nu_1}})*\dots * \globP(\mathbf{D}^{n_{\nu_n}}) \longrightarrow X_\lambda
\]
by ``concatenating'' the attaching maps of the globular cells $c_{\nu_1},\dots,c_{\nu_n}$. Let $\alpha_{i-1}$ ($\alpha_{i}$ resp.) be the initial state (the final state resp.) of $\globP(\mathbf{D}^{n_{\nu_i}})$ for $1\leqslant i\leqslant n$ in $\globP(\mathbf{D}^{n_{\nu_1}})* \dots * \globP(\mathbf{D}^{n_{\nu_n}})$. It induces a continuous map 
\[
\PA^{top}\widehat{g_{\underline{c}}}:X_{\underline{c}}=\PA^{top}_{\alpha_0,\alpha_n}(\globP(\mathbf{D}^{n_{\nu_1}})*\dots * \globP(\mathbf{D}^{n_{\nu_n}})) \longrightarrow \PA^{top} X_\lambda.
\]
\end{nota}

As a consequence of the associativity of the semimonoidal structure on $\mathcal{P}$-spaces and of \cite[Proposition~5.16]{Moore1}, we have

\bp \label{Ftenseur} 
Let $U_1,\dots,U_p$ be $p$ topological spaces with $p\geq 1$. Let $\ell_1,\dots,\ell_p>0$. There is the natural isomorphism of $\mathcal{P}$-spaces 
\[
\mathbb{F}^{\mathcal{P}^{op}}_{\ell_1}U_1\ot \dots \ot \mathbb{F}^{\mathcal{P}^{op}}_{\ell_p}U_p \iso \mathbb{F}^{\mathcal{P}^{op}}_{\ell_1+\dots+\ell_p}(U_1\p \dots\p U_p).
\]
\ep

\bp  \label{comp-gl} (\cite[Proposition~6.3]{Moore2} for $\mathcal{G}$ and $\mathcal{M}$)
Let $Z_1,\dots,Z_p$ be $p$ topological spaces of $\topspace$ with $p\geq 1$. Consider the $\mathcal{P}$-multipointed $d$-space $X=\globP(Z_1)* \dots *\globP(Z_p)$ with $p\geq 1$. There is a homeomorphism 
\[
\PA_{\alpha_0,\alpha_p}^{top} X \iso  \mathcal{P}(1,p)\p Z_1\p \dots \p Z_p.
\]
\ep

The case $p=1$ is treated in Proposition~\ref{calcul-topology-glob}. The proof of Proposition~\ref{comp-gl} is a modified version of the proof of \cite[Proposition~6.3]{Moore2}, the latter working only for the case $\mathcal{P}=\mathcal{G}$. Like in the proof of Proposition~\ref{calcul-topology-glob}, the verification of the continuity in one direction is different from the proof of \cite[Proposition~6.3]{Moore2} because of the possible presence, in the case $\mathcal{P}=\mathcal{M}$, of stop intervals.

\bpf
The Moore composition of paths induced a map of $\mathcal{P}$-spaces
\[
\PA_{0,1}^\bullet\globP(Z_1) \ot \dots \ot \PA_{0,1}^\bullet\globP(Z_p) \longrightarrow \PA_{\alpha_0,\alpha_p}^\bullet X.
\]
By Proposition~\ref{calcul-topology-glob}, there is the isomorphism of $\mathcal{P}$-spaces 
\[
\PA_{0,1}^\bullet\globP(Z) \iso \mathbb{F}_1^{\mathcal{P}^{op}}Z
\]
for all topological spaces $Z$. We obtain a map of $\mathcal{P}$-spaces
\[
\mathbb{F}_1^{\mathcal{P}^{op}}Z_1 \ot \dots \ot \mathbb{F}_1^{\mathcal{P}^{op}}Z_p \longrightarrow \PA_{\alpha_0,\alpha_p}^\bullet X.
\]
By Proposition~\ref{Ftenseur}, and since $\PA_{\alpha_0,\alpha_p}^1 X=\PA_{\alpha_0,\alpha_p}^{top} X$ by definition of the functor $\PA_{\alpha_0,\alpha_p}^\bullet X$, we obtain a continuous map 
\[
\begin{cases}
	\Psi:&\mathcal{P}(1,p)\p Z_1\p \dots \p Z_p  \longrightarrow \PA_{\alpha_0,\alpha_p}^{top} X\\
	&(\phi,z_1,\dots,z_p) \mapsto (\delta_{z_1}\phi_1)*\dots *(\delta_{z_p}\phi_p)
\end{cases}
\]
where $\phi_i\in \mathcal{P}(\ell_i,1)$ with $\sum_i \ell_i=1$ and $\phi=\phi_1\ot \dots \ot \phi_p$ being a decomposition given by the third axiom of reparametrization category. The map $\Psi$ is bijective by Theorem~\ref{normal-form}. The continuous maps $Z_i\to \{0\}$ for $1\leqslant i \leqslant p$ induce by functoriality a map of $\mathcal{P}$-multipointed $d$-spaces $X \to \vI^{\mathcal{P}}*\dots * \vI^{\mathcal{P}}$ ($p$ times) and then a continuous map \[
\begin{cases}
	k:&\PA_{\alpha_0,\alpha_p}^{top}X\longrightarrow \PA_{\alpha_0,\alpha_p}^{top}(\vI^{\mathcal{P}}*\dots * \vI^{\mathcal{P}}) = \mathcal{P}(1,p)\\
	&(\delta_{z_1}\phi_1)*\dots *(\delta_{z_p}\phi_p)\mapsto (\delta_{0}\phi_1)*\dots *(\delta_{0}\phi_p) = \phi_1\ot \dots \ot \phi_p.
\end{cases}\] 
Consider the set map 
\[
\begin{cases}
	\overline{k}:&\PA_{\alpha_0,\alpha_p}^{top}X\longrightarrow Z_1\p  \dots \p Z_p\\
	& (\delta_{z_1}\phi_1)*\dots *(\delta_{z_p}\phi_p) \mapsto (z_1,\dots,z_p).
\end{cases}
\]
Let $i\in \{1,\dots,p\}$. Suppose that the composite set map $\pr_i\overline{k}:\PA_{\alpha_0,\alpha_p}^{top} \to Z_i$ is not continuous where $\pr_i$ is the projection on the $i$-th factor. All involved topological spaces being sequential, there exist $z_i^\infty\in Z_i$, an open neighborhood $V$ of $z_i^\infty$ in $Z_i$, and a sequence $((\delta_{z^n_1}\phi^n_1)*\dots *(\delta_{z^n_p}\phi^n_p))_{n\geq 0}$ which converges to $(\delta_{z^\infty_1}\phi^\infty_1)*\dots *(\delta_{z^\infty_p}\phi^\infty_p)$ such that $z^n_i\in Z_i\backslash V$ for all $n\geq 0$. Let $\phi^n=\phi_1^n \ot \dots \ot \phi_p^n$ for $n\geq 0$ and $\phi^\infty=\phi_1^\infty \ot \dots \ot \phi_p^\infty$. Choose $t_0\in [0,1]$ such that $\phi^\infty(t_0)\in ]i-1,i[$. The sequence \[\big(\big((\delta_{z^n_1}\phi^n_1)*\dots *(\delta_{z^n_p}\phi^n_p)\big)(t_0)\big)_{n\geq 0}\] converges to \[((\delta_{z^\infty_1}\phi^\infty_1)*\dots *(\delta_{z^\infty_p}\phi^\infty_p))(t_0) = (\delta_{z^\infty_1}*\dots *\delta_{z^\infty_p})(\phi^\infty(t_0)).\] By continuity of the map $k:\PA_{\alpha_0,\alpha_p}^{top}X\to \mathcal{P}(1,p)$, the sequence $(\phi^n(t_0))_{n\geq 0}$ of $[0,p]$ converges to $\phi^\infty(t_0) \in ]i-1,i[$. This implies that there exists $N\geq 0$ such that for all $n\geq N$, $\phi^n(t_0)\in ]i-1,i[$. We obtain that the sequence $((z_i^n,\phi^n(t_0)-i+1))_{n\geq N}$ converges to $(z_i^\infty,\phi^\infty(t_0)-i+1)$ in $|\globP(Z_i)|\backslash \{0,1\}$. By considering the well-defined projection (the left-hand term being equipped with the relative topology) \[|\globP(Z_i)|\backslash \{0,1\}\longrightarrow Z_i,\] we obtain that the sequence $(z_i^n)_{n\geq N}$ converges to $z_i^\infty$, and therefore that $z_i^\infty\in Z_i\backslash V$, the latter set being closed in $Z_i$: contradiction. 

This means that the composite set map $\pr_i\overline{k}$ is continuous for all $i\in \{1,\dots,p\}$, and therefore, by the universal property of the product, that the set map $\overline{k}:\PA_{\alpha_0,\alpha_p}^{top}X\to Z_1\p  \dots \p Z_p$ is continuous. This implies that the set map 
\[
\Psi^{-1}=(k,\overline{k}):(\delta_{z_1}\phi_1)*\dots *(\delta_{z_p}\phi_p) \mapsto (\phi_1\ot \dots \ot \phi_p,z_1,\dots,z_p).
\]
is continuous and that $\Psi$ is a homeomorphism. 
\epf

\begin{cor}
	Let $Z_1,\dots,Z_p$ be $p$ topological spaces of $\topspace$ with $p\geq 1$. Consider the $\mathcal{P}$-multipointed $d$-space $X=\globP(Z_1)* \dots *\globP(Z_p)$ with $p\geq 1$. There is a homeomorphism 
	\[
	\PA_{\alpha_0,\alpha_p} X \iso  Z_1\p \dots \p Z_p.
	\]
\end{cor}

\bpf
There are the homeomorphisms 
\[
\liminj_{\mathcal{P}^1} \bigg(\mathcal{P}(1,p)\p Z_1\p \dots \p Z_p\bigg) \iso \bigg(\liminj_{\mathcal{P}^1} \mathcal{P}(1,p)\bigg)\p Z_1\p \dots \p Z_p \iso Z_1\p \dots \p Z_p,
\]
the left-hand homeomorphism since $\topspace$ is cartesian closed and the right-hand homeomorphism by Proposition~\ref{pre-un}. The proof is complete thanks to Proposition~\ref{comp-gl} and Proposition~\ref{un}.
\epf

\begin{lem} \label{prod-seq} \cite[Lemma~6.10]{Moore2}
	Let $U_1,\dots,U_p$ be $p$ first-countable $\Delta$-Hausdorff $\Delta$-gene\-rated spaces with $p\geq 1$. Let $(u_n^i)_{n\geq 0}$ be a sequence of $U_i$ for $1\leqslant i \leqslant p$ which converges to $u_\infty^i\in U_i$. Then the sequence $((u_n^1,\dots,u_n^p))_{n\geq 0}$ converges to $(u_\infty^1,\dots,u_\infty^p)\in U_1\p \dots \p U_p$ for the product calculated in $\topspace$.
\end{lem}

\begin{nota}
	Let $\underline{c}$ be the carrier of some execution path of $X_\lambda$. Using the identification provided by the homeomorphism of Proposition~\ref{comp-gl}, we can use the notation
	\[
	(\PA^{top}\widehat{g_{\underline{c}}})(\phi,z^1,\dots,z^n) = (\widehat{g_{\nu_1}}\delta_{z^{1}} * \dots * \widehat{g_{\nu_n}}\delta_{z^{n}})\phi. 
	\]
\end{nota}

Lemma~\ref{converge} is implicitly used in \cite[Theorem~6.11 and Theorem~7.3]{Moore2} and in \cite[Theorem~7.7]{NaturalRealization}.

\begin{lem} \label{converge}
	Let $X$ be a sequential topological space. Let $x_\infty\in X$. Let $(x_n)_{n\geq 0}$ be a sequence such that $x_\infty$ is a limit point of all subsequences. Then the sequence $(x_n)_{n\geq 0}$ converges to $x_\infty$.
\end{lem}

\bpf
Otherwise, consider an open neighborhood $V$ of $x_\infty$ such that for all $n\geq N$, $x_n\notin V$ for some $N\geq 0$. This means that $x_\infty\in X\backslash V$, the subset $X\backslash V$ being sequentially closed in $X$: contradiction.
\epf

\bth \label{img-closed} Let $\underline{c}$ be the carrier of some execution path of $X_\lambda$. 
\begin{enumerate}[leftmargin=*]
	\item Consider a sequence $(\gamma_k)_{k\geq 0}$ in the image of $\PA^{top}\widehat{g_{\underline{c}}}$ which converges pointwise to $\gamma_\infty$ in $\PA^{top}X_\lambda$. Let \[\gamma_k=(\PA^{top}\widehat{g_{\underline{c}}})(\phi_k,z_k^1,\dots,z_k^n)\]
	with $\phi_k\in \mathcal{P}(1,n)$ and $z_k^i\in \mathbf{D}^{n_{\nu_i}}$ for $1\leqslant i \leqslant n$ and $k\geq 0$. Then there exist $\phi_\infty\in \mathcal{P}(1,n)$ and $z_\infty^i\in \mathbf{D}^{n_{\nu_i}}$ for $1\leqslant i \leqslant n$ such that 
	\[\gamma_\infty=(\PA^{top}\widehat{g_{\underline{c}}})(\phi_\infty,z_\infty^1,\dots,z_\infty^n)\]
	and such that $(\phi_\infty,z_\infty^1,\dots,z_\infty^n)$ is a limit point of the sequence $((\phi_k,z_k^1,\dots,z_k^n))_{k\geq 0}$.
	\item The image of $\PA^{top}\widehat{g_{\underline{c}}}$ is closed in $\PA^{top} X_\lambda$.
\end{enumerate}
\eth

\bpf 
The case $\mathcal{P}=\mathcal{G}$ is treated in \cite[Theorem~6.11]{Moore2}. Let us suppose that $\mathcal{P}=\mathcal{M}$. The proof is similar but simpler because it is not necessary to verify anymore that some limit execution paths are regular. 

(1) By a Cantor diagonalization argument, we can suppose that the sequence $(z_k^i)_{k\geq 0}$ converges to $z_\infty^i\in \mathbf{D}^{n_{\nu_i}}$ for each $1\leqslant i \leqslant n$ and that the sequence $(\phi_k(r))_{k\geq 0}$ converges to a real number denoted by $\phi_\infty(r)\in [0,n]$ for all rational numbers $r\in [0,1]\cap \mathbb{Q}$. Since the sequence of execution paths $(\gamma_k)_{k\geq 0}$ converges pointwise to $\gamma_\infty$, we obtain 
\[
\gamma_\infty(r) = (\widehat{g_{\nu_1}}\delta_{z_\infty^{1}} * \dots * \widehat{g_{\nu_n}}\delta_{z_\infty^{n}})(\phi_\infty(r))
\]
for all $r\in [0,1]\cap \mathbb{Q}$. For $r_1<r_2 \in [0,1]\cap \mathbb{Q}$, $\phi_k(r_1)\leqslant \phi_k(r_2)$ for all $k\geq 0$. Therefore by passing to the limit, we obtain $\phi_\infty(r_1)\leqslant \phi_\infty(r_2)$. Note that $\phi_\infty(0)=0$ and $\phi_\infty(1)=n$ since $0,1\in \mathbb{Q}$. For $t\in ]0,1[$, let us extend the definition of $\phi_\infty$ as follows: 
\[
\phi_\infty(t) = \sup \{\phi_\infty(r)\mid r\in ]0,t]\cap \mathbb{Q}\}.
\]
By continuity, we deduce that \[
\gamma_\infty(t) = (\widehat{g_{\nu_1}}\delta_{z_\infty^{1}} * \dots * \widehat{g_{\nu_n}}\delta_{z_\infty^{n}})(\phi_\infty(t))
\]
for all $t\in [0,1]$. It is easy to see that the set map $\phi_\infty:[0,1]\to [0,n]$ is nondecreasing and that it preserves extremities. By definition of the Moore composition, there exist $0=t_0\leqslant t_1 \leqslant \dots \leqslant t_n=1$ such that for all $1\leqslant i\leqslant n$,
\[
\forall t\in [t_{i-1},t_i], \gamma_\infty(t) = \widehat{g_{\nu_i}}(z_\infty^{i},\phi_\infty(t)-i+1).
\]
This implies that the restriction of $\phi_\infty$ to $[t_{i-1},t_i]$ is surjective. We deduce that the nondecreasing set map $\phi_\infty:[0,1] \to [0,n]$ is surjective, and therefore that $\phi_\infty\in \mathcal{M}(1,n)$. Let $t\in [0,1]\backslash \mathbb{Q}$. The sequence $(\phi_k(t))_{k\geq 0}$ has at least one limit point $\ell$. There exists a subsequence of $(\phi_k(t))_{k\geq 0}$ which converges to $\ell$. We obtain: $\forall r\in [0,t]\cap \mathbb{Q},\forall r'\in [t,1]\cap \mathbb{Q}, \phi_\infty(r) \leqslant \ell\leqslant \phi_\infty(r')$. Since $\phi_\infty\in \mathcal{M}(1,n)$ and by density of $\mathbb{Q}$, we deduce that $\ell=\phi_\infty(t)$ necessarily. Using Lemma~\ref{converge}, we deduce that the sequence $(\phi_k)_{k\geq 0}$ converges pointwise to $\phi_\infty$. Using Proposition~\ref{morphG-metrizable}, we deduce that $(\phi_k)_{k\geq 0}$ converges uniformly to $\phi_\infty$. We deduce that $(\phi_\infty,z_\infty^1,\dots,z_\infty^n)$ is a limit point of the sequence $((\phi_k,z_k^1,\dots,z_k^n))_{k\geq 0}$ in $\mathcal{M}(1,n)\p \mathbf{D}^{n_{\nu_1}}\p \dots \p \mathbf{D}^{n_{\nu_n}}$ by Proposition~\ref{morphG-metrizable} and Lemma~\ref{prod-seq}.

(2) Let $(\PA^{top}\widehat{g_{\underline{c}}}(\Gamma_n))_{n\geq 0}$ be a sequence of $(\PA^{top}\widehat{g_{\underline{c}}})(X_{\underline{c}})$ which converges in $\PA^{top} X_\lambda$. The limit $\gamma_\infty\in \PA^{top} X_\lambda$ of the sequence of execution paths $(\PA^{top}\widehat{g_{\underline{c}}}(\Gamma_n))_{n\geq 0}$ is also a pointwise limit. We can suppose by extracting a subsequence that the sequence $(\Gamma_n)_{n\geq 0}$ of $X_{\underline{c}}$ converges in $X_{\underline{c}}$. Thus, by continuity of $\PA^{top}\widehat{g_{\underline{c}}}$, we obtain $\gamma_\infty=(\PA^{top}\widehat{g_{\underline{c}}})(\Gamma_\infty)$ for some $\Gamma_\infty\in X_{\underline{c}}$. We deduce that $\PA^{top}\widehat{g_{\underline{c}}}(X_{\underline{c}})$ is sequentially closed in $\PA^{top} X_\lambda$. Since $\PA^{top} X_\lambda$ is sequential, being a $\Delta$-generated space, the proof is complete. 
\epf

As a corollary of Theorem~\ref{img-closed}, we obtain:

\begin{cor} \label{Dini-cellular}
	Suppose that $X_\lambda$ is a finite cellular $\mathcal{P}$-multipointed  $d$-space, i.e. $X^0$ is finite and $\lambda$ is a finite ordinal. If $X_\lambda$ has no loops, then the topology of $\PA^{top}X_\lambda$ is the topology of the pointwise convergence which is therefore $\Delta$-generated.
\end{cor}

\bpf
It is mutatis mutandis the proof of \cite[Corollary~6.12]{Moore2}.
\epf

\section{Locally finite cellular multipointed d-space}
\label{locallyfinite}

We want to give an application of Theorem~\ref{diagonal-execution} and Theorem~\ref{img-closed} before addressing the main subject of this paper. The reading of this section is not necessary to understand Section~\ref{unit}. A cellular $\mathcal{P}$-multipointed $d$-space $X_\lambda$ is fixed. 

\bd \label{lf}
The cellular $\mathcal{P}$-multipointed $d$-space space $X_\lambda$ is \textit{locally finite} if for all $\nu<\lambda$, the set $\{\nu'<\lambda\mid  \widehat{c_{\nu'}}\cap c_\nu \neq \varnothing\}$ is finite and each state meets a finite number of $\widehat{c_{\nu}}$. In other terms, the underlying topological space $|X_\lambda|$, which is cellular by Proposition~\ref{p1}, is locally finite.
\ed

Lemma~\ref{lpc} is a consequence of \cite[Proposition~3.4]{MR3270173} and \cite[Proposition~3.10]{MR3270173}. It can be easily proved without using diffeological spaces.

\begin{lem} (well-known) \label{lpc}
	Every $\Delta$-generated space is locally path-connected. 
\end{lem}

\bpf
Let $U$ be an open subset of a $\Delta$-generated space $X$. Then $U$ equipped with the relative topology is $\Delta$-generated by \cite[Proposition~2.4]{leftproperflow}. Therefore $U$ equipped with the relative topology is homeomorphic to the disjoint sum of its path-connected components by \cite[Proposition~2.8]{mdtop}. Thus $X$ is locally path-connected.
\epf

\bd 
A topological space $X$ is \textit{weakly locally path-connected} if for every $x\in X$ and every neighborhood $W$ of $x$, there exists a path-connected neighborhood (not necessarily open) $W'$ of $x$ such that $W'\subset W$.
\ed

\begin{lem} (well-known) \label{lpathc}
	Every weakly locally path-connected space is locally path-connec\-ted.
\end{lem}

\bpf
Let $W$ be a neighborhood of $x\in X$. Then there exists a path-connected neighborhood $W'$ of $x$ such that $W'\subset W$. This means that $W'$ is included in the path-connected component $C$ of $x$ in $W$. Therefore $x$ is in the interior of $C$. Thus $C$ is open and $X$ is locally path-connected.
\epf

\bp \label{met}
Let $\lambda$ be an ordinal. Let $Z:\lambda\to \topspace$ be a colimit-preserving functor such that $Z_\lambda$ is cellular for the q-model structure of $\topspace$. If the cellular space $Z_\lambda$ is locally finite, then the topological space $Z_\lambda$ is metrizable.
\ep

\bpf[Sketch of proof]
The technique used in \cite{MR1074175} to reorganize and regroup the cells in a CW-complex using the notion of star of a subset \cite[Example~2]{MR1074175} works in the same way for cellular topological spaces, even when $\lambda$ is not countable. Assume first that $Z_\lambda$ is path-connected. By \cite[Proposition~1.5.12]{MR1074175}, the ordinal $\lambda$ is countable, $Z_\lambda$ being locally finite. Using \cite[Proposition~1.5.13]{MR1074175}, the cells are reorganized so that $\lambda=\aleph_0$ and so that each $Z_n$ for $n$ finite is a finite cellular topological space (i.e. built using finitely many cells). Moreover, for all $n\geq 0$, the space $Z_n$ is contained in the interior $\overset{\circ}{Z_{n+1}}$ of $Z_{n+1}$ for the topology of $Z_\lambda$ and there is the equality \[\bigcup_{n\geq 0}\overset{\circ}{Z_n} = \bigcup_{n\geq 0}Z_n= Z_\lambda.\] 
Using \cite[Theorem~1.5.16]{MR1074175}, we deduce that $Z_\lambda$ is metrizable and that it can be embedded in the Hilbert cube equipped with the $\ell^2$ metric. This means that the metric of $Z_\lambda$ is bounded, namely by the constant $\pi/\sqrt{6}$ which does not depend on $Z_\lambda$. In the general case, the $\Delta$-generated space $Z_\lambda$ is homeomorphic to the disjoint sum of its path-connected components by \cite[Proposition~2.8]{mdtop}. Thus, the metric on each path-connected component being bounded by $\pi/\sqrt{6}$, the disjoint sum is metrizable.
\epf

In fact, we could prove the equivalence for cellular topological spaces of the conditions locally finite, metrizable, locally compact and first-countable as it is done in \cite[Proposition~1.5.10 and Proposition~1.5.17]{MR1074175} for CW-complexes.

\begin{cor} \label{met2}
	Assume $X_\lambda$ locally finite. Let $(\alpha,\beta)\in X^0\p X^0$. Then the topological space $(\PA^{top}_{\alpha,\beta}X_\lambda)_{co}$ is metrizable, and therefore sequential and first-countable.
\end{cor}

\bpf
By Proposition~\ref{met}, the topological space $|X_\lambda|$ is metrizable, $X_\lambda$ being locally finite by hypothesis. The space $(\PA^{top}_{\alpha,\beta}X_\lambda)_{co}$ is therefore metrizable by \cite[Proposition~A.13]{MR1867354}.
\epf

\bp \label{locally-path-connected}
Assume $X_\lambda$ locally finite. Let $(\alpha,\beta)\in X^0\p X^0$. The space $(\PA^{top}_{\alpha,\beta}X_\lambda)_{co}$ is locally path-connected.
\ep

\bpf 
By Lemma~\ref{lpathc}, it suffices to prove that $(\PA^{top}_{\alpha,\beta}X_\lambda)_{co}$ is weakly locally path-connected. Consider an execution path $\gamma_\infty$ of $\PA^{top}_{\alpha,\beta}X_\lambda$. Let $\Omega$ be an open subset of $(\PA^{top}_{\alpha,\beta}X_\lambda)_{co}$ containing $\gamma_\infty$. Then $\Omega$ is an open subset of $\PA^{top}_{\alpha,\beta}X_\lambda$ since the $\Delta$-kelleyfication adds open subsets. Let $\mathcal{T}$ be the set of all carriers of all execution paths of $\PA^{top}_{\alpha,\beta}X_\lambda$. Consequently, for each carrier $\underline{c} \in \mathcal{T}$ and for each $\Gamma\in (\PA^{top} \widehat{g_{\underline{c}}})^{-1}(\gamma_\infty)$, there exists an open neighborhood $\Omega_\Gamma$ of $\Gamma$ such that $(\PA^{top} \widehat{g_{\underline{c}}})(\Omega_\Gamma) \subset \Omega$. By Lemma~\ref{lpc}, we can suppose that $\Omega_\Gamma$ is path-connected. Consider 
\[
U = \bigcup_{\underline{c}\in \mathcal{T}} \bigcup_{\Gamma\in (\PA^{top} \widehat{g_{\underline{c}}})^{-1}(\gamma_\infty)} (\PA^{top} \widehat{g_{\underline{c}}})(\Omega_\Gamma).
\]
Then $U$ is path-connected and $U\subset \Omega$. Suppose that $\gamma_\infty$ is not in the interior of $U$ in $(\PA^{top}_{\alpha,\beta}X_\lambda)_{co}$. The space $(\PA^{top}_{\alpha,\beta}X_\lambda)_{co}$ being sequential by Corollary~\ref{met2}, there exists a sequence $(\gamma_n)_{n\geq 0}$ of execution paths not belonging to $U$  converging to $\gamma_\infty$ in $(\PA^{top}_{\alpha,\beta}X_\lambda)_{co}$. Since the set $\{\carrier(\gamma_n)\mid n\geq 0\}$ is finite by Theorem~\ref{diagonal-execution}, we can always suppose that the sequence of carriers $(\carrier(\gamma_n))_{n\geq 0}$ is constant and e.g. equal to some $\underline{c}\in \mathcal{T}$ by extracting a subsequence. Therefore we can write $\gamma_n = (\PA^{top} \widehat{g_{\underline{c}}})(\Gamma_n)$ with $\Gamma_n\in X_{\underline{c}}$ (see Notation~\ref{Xc}). The sequence of execution paths $(\gamma_n)_{n\geq 0}$ converges pointwise to $\gamma_\infty$. Thus, by Theorem~\ref{img-closed}, we can suppose that the sequence $(\Gamma_n)_{n\geq 0}$ converges to $\Gamma_\infty\in X_{\underline{c}}$ after extracting a subsequence again. By continuity, we obtain the equality $\gamma_\infty = (\PA^{top} \widehat{g_{\underline{c}}})(\Gamma_\infty)$. There exists $N\geq 0$ such that for all $n\geq N$, $\Gamma_n\in \Omega_{\Gamma_\infty}$, i.e. $\gamma_n = (\PA^{top} \widehat{g_{\underline{c}}})(\Gamma_n) \in U$ for all $n\geq N$. Contradiction. Thus $\gamma_\infty$ is in the interior of $U$. 
\epf

\bth \label{main} Assume $X_\lambda$ locally finite. Let $(\alpha,\beta)\in X^0\p X^0$. The topological space $(\PA^{top}_{\alpha,\beta}X_\lambda)_{co}$ equipped with the compact-open topology is $\Delta$-generated. The topological space \[\PA^{top}X_\lambda = \coprod_{(\alpha,\beta)\in X^0\p X^0} \PA^{top}_{\alpha,\beta}X_\lambda\] is metrizable with the distance of the uniform convergence. The underlying topology is the compact-open topology.
\eth

\bpf 
By Corollary~\ref{met2}, the topological space $(\PA^{top}_{\alpha,\beta}X_\lambda)_{co}$ is first-countable. The space $(\PA^{top}_{\alpha,\beta}X_\lambda)_{co}$ is locally path-connected by Proposition~\ref{locally-path-connected}. Using \cite[Proposition~3.11]{MR3270173}, we deduce that $(\PA^{top}_{\alpha,\beta}X_\lambda)_{co}$ is $\Delta$-generated. The set of all execution paths equipped with the compact-open topology $(\PA^{top}X_\lambda)_{co}$ satisfies 
\[
(\PA^{top}X_\lambda)_{co} \iso \coprod_{(\alpha,\beta)\in X^0\p X^0} (\PA^{top}_{\alpha,\beta}X_\lambda)_{co}
\]
because $X^0$ is a discrete subspace of $|X_\lambda|$. Hence $(\PA^{top}X_\lambda)_{co}$ is $\Delta$-generated  and metrizable by the distance of the uniform convergence.
\epf

\section{Multipointed d-space and Moore flow}
\label{unit}

Consider a pushout diagram of $\mathcal{P}$-multipointed $d$-spaces 
\[
\begin{tikzcd}[row sep=3em, column sep=3em]
	\globP(\mathbf{S}^{n-1}) \arrow[d] \arrow[r,"g"] & A \arrow[d,"f"] \\
	\globP(\mathbf{D}^{n}) \arrow[r,"\widehat{g}"] & \cocartesian X
\end{tikzcd}
\]
with $n\geq 0$ and $A$ cellular. Note that $A^0=X^0$. Let $D= \mathbb{F}^{\mathcal{P}^{op}}_1\mathbf{S}^{n-1}$ and $E=\mathbb{F}^{\mathcal{P}^{op}}_1\mathbf{D}^{n}$. Consider the $\mathcal{P}$-flow $\overline{X}$ defined by the pushout diagram of Figure~\ref{Xoverline} where the two equalities 
\[
\begin{aligned}
	&\moore^{\mathcal{P}}(\globP(\mathbf{S}^{n-1})) = \glob(D)\\
	&\moore^{\mathcal{P}}(\globP(\mathbf{D}^{n})) = \glob(E)
\end{aligned}
\]
come from Proposition~\ref{calculM} and where the map $\psi$ is induced by the universal property of the pushout. 

\begin{figure}
	\[
	\begin{tikzcd}[row sep=2em, column sep=3em]
		\moore^{\mathcal{P}}(\globP(\mathbf{S}^{n-1})) = \glob(D) \arrow[d] \arrow[r,"\moore^{\mathcal{P}}(g)"] & \moore^{\mathcal{P}}(A) \arrow[rdd,bend left=20pt,"\moore^{\mathcal{P}}(f)"] \arrow[d,"\overline{f}"] \\
		\moore^{\mathcal{P}}(\globP(\mathbf{D}^{n})) = \glob(E) \arrow[r,"\overline{g}"] \arrow[rrd,bend right=20pt,"\moore^{\mathcal{P}}(\widehat{g})"] & \cocartesian \overline{X}\arrow[rd,dashed,"\psi"]&  \\
		&& \moore^{\mathcal{P}}(X).
    \end{tikzcd}
	\]
	\caption{Definition of $\overline{X}$}
	\label{Xoverline}
\end{figure}

The $\mathcal{P}$-space of execution paths of the Moore flow $\overline{X}$ can be calculated by introducing a diagram of $\mathcal{P}$-spaces $\mathcal{D}^{f}$ over a Reedy category $\mathcal{P}^{g(0),g(1)}(A^0)$ whose definition is recalled now. It was introduced for the first time in \cite[Section~3]{leftproperflow}. 

Let $S$ be a nonempty set. Let $\mathcal{P}^{u,v}(S)$~\footnote{The use of the letter $\mathcal{P}$ here has nothing to do with the reparametrization category $\mathcal{P}$. It is a bit unfortunate but I prefer not to change the notation.} be the small category defined by generators and relations as follows: 
\begin{itemize}
	\item $u,v\in S$ ($u$ and $v$ may be equal).
	\item The objects are the tuples of the form 
	\[\underline{m}=((u_0,\epsilon_1,u_1),(u_1,\epsilon_2,u_2),\dots ,(u_{n-1},\epsilon_n,u_n))\]
	with $n\geq 1$, $u_0,\dots,u_n \in S$, $\epsilon_1,\dots,\epsilon_n \in \{0,1\}$ and \[\forall i\hbox{ such that } 1\leqslant i\leqslant n, \epsilon_i = 1\Rightarrow (u_{i-1},u_i)=(u,v).\] 
	\item There is an arrow \[c_{n+1}:(\underline{m},(x,0,y),(y,0,z),\underline{n}) \to (\underline{m},(x,0,z),\underline{n})\]
	for every tuple $\underline{m}=((u_0,\epsilon_1,u_1),(u_1,\epsilon_2,u_2),\dots ,(u_{n-1},\epsilon_n,u_n))$ with $n\geq 1$ and every tuple $\underline{n}=((u'_0,\epsilon'_1,u'_1),(u'_1,\epsilon'_2,u'_2),\dots ,(u'_{n'-1},\epsilon'_{n'},u'_{n'}))$ with $n'\geq 1$. It is called a \textit{composition map}. 
	\item There is an arrow \[I_{n+1}:(\underline{m},(u,0,v),\underline{n}) \to (\underline{m},(u,1,v),\underline{n})\] for every tuple $\underline{m}=((u_0,\epsilon_1,u_1),(u_1,\epsilon_2,u_2),\dots ,(u_{n-1},\epsilon_n,u_n))$ with $n\geq 1$ and every tuple $\underline{n}=((u'_0,\epsilon'_1,u'_1),(u'_1,\epsilon'_2,u'_2),\dots ,(u'_{n'-1},\epsilon'_{n'},u'_{n'}))$ with $n'\geq 1$.
	It is called an \textit{inclusion map}. 
	\item There are the relations (group A) $c_i.c_j = c_{j-1}.c_i$ if $i<j$ (which means since $c_i$ and $c_j$ may correspond to several maps that if $c_i$ and $c_j$ are composable, then there exist $c_{j-1}$ and $c_i$ composable satisfying the equality). 
	\item There are the relations (group B) $I_i.I_j = I_j.I_i$ if $i\neq j$. By definition of these maps, $I_i$ is never composable with itself. 
	\item There are the relations (group C) \[c_i.I_j = \begin{cases}
		I_{j-1}.c_i&\hbox{if } j\geq i+2\\
		I_j.c_i&\hbox{if } j\leqslant i-1.
	\end{cases}\]
	By definition of these maps, $c_i$ and $I_i$ are never composable as well as $c_i$ and $I_{i+1}$. 
\end{itemize}
By \cite[Proposition~3.7]{leftproperflow}, there exists a structure of Reedy category on $\mathcal{P}^{u,v}(S)$ with the $\mathbb{N}$-valued degree map defined by \[d((u_0,\epsilon_1,u_1),(u_1,\epsilon_2,u_2),\dots ,(u_{n-1},\epsilon_n,u_n)) = n + \sum_i \epsilon_i.\]
The maps raising the degree are the inclusion maps. The maps decreasing the degree are the composition maps.

Let $T$ be the $\mathcal{P}$-space defined by the pushout diagram of $\topdgrP_0$
\[
\begin{tikzcd}[row sep=3em, column sep=4em]
		D  \arrow[d] \arrow[r,"\PA\moore^{\mathcal{P}}(g)"] & \PA_{g(0),g(1)} \moore^{\mathcal{P}}(A) \arrow[d,"\PA\overline{f}"] \\
	E  \arrow[r,"\PA\overline{g}"] & \cocartesian T.
\end{tikzcd}
\]
Consider the diagram of spaces $\D^{f}:\mathcal{P}^{g(0),g(1)}(A^0)\to \topdgrP_0$ defined as follows:
\[
\D^f((u_0,\epsilon_1,u_1),(u_1,\epsilon_2,u_2),\dots ,(u_{n-1},\epsilon_n,u_n)) = Z_{u_0,u_1}\ot Z_{u_1,u_2} \ot \dots \ot Z_{u_{n-1},u_n}
\]
with 
\[
Z_{u_{i-1},u_i}=
\begin{cases}
	\PA_{u_{i-1},u_i}\moore^{\mathcal{P}}(A) & \hbox{if }\epsilon_i=0\\
	T & \hbox{if }\epsilon_i=1
\end{cases}
\] 
In the case $\epsilon_i=1$, $(u_{i-1},u_i)=(g(0),g(1))$ by definition of $\mathcal{P}^{g(0),g(1)}(A^0)$. The inclusion maps $I_i's$ are induced by the map $\PA\overline{f}:\PA_{g(0),g(1)} \moore^{\mathcal{P}}(A) \to T$. The composition maps $c_i's$ are induced by the compositions of paths of the Moore flow $\moore^{\mathcal{P}}(A)$.

\bth \cite[Theorem~9.7]{Moore1} \label{calcul_pushout}
We obtain a well-defined diagram of $\mathcal{P}$-spaces \[\D^f:\mathcal{P}^{g(0),g(1)}(A^0)\to \topdgrP_0.\] There is the isomorphism of $\mathcal{P}$-spaces $\liminj \D^f \iso \PA \overline{X}$. 
\eth

By the universal property of the pushout, we obtain a canonical map of $\mathcal{P}$-spaces \[\PA \psi:\liminj \D^f \longrightarrow \PA\mathbb{M}^\mathcal{P}X.\]

\bd 
Let $\underline{x}$ be an element of some vertex of the diagram of spaces $\mathcal{D}^f$. We say that $\underline{x} \in \mathcal{D}^f (\underline{n})$ is \textit{simplified} if
\[
d(\underline{n}) = \min \big\{d(\underline{m}) \mid \exists \underline{m}\in \Obj(\mathcal{P}^{g({0}),g({1})}(A^0)) \hbox{ and }\exists \underline{y}\in \mathcal{D}^f (\underline{m}),\underline{y}=\underline{x} \in \liminj \mathcal{D}^f 
\big\}.
\]
\ed

\bth \label{pre-calculation-pathspace}
Under the hypotheses and the notations of this section. The map of $\mathcal{P}$-spaces \[\PA\psi:\liminj \D^f \longrightarrow \PA\mathbb{M}^\mathcal{P}(X)\] is an isomorphism.
\eth

\bpf
The structure of the proof is the same as the one of the proofs of \cite[Theorem~7.2 and Theorem~7.3]{Moore2}. At first it must be proved that the map $\PA\psi$ is an objectwise bijection. The role of \cite[Theorem~5.20]{Moore2} is played by Theorem~\ref{calcul_final_structure}. Then it must be proved that the map $\PA\psi$ is an objectwise homeomorphism. The roles of \cite[Theorem~5.19]{Moore2} and \cite[Theorem~6.11]{Moore2} are played by Theorem~\ref{carrier-finite-on-compact} and Theorem~\ref{img-closed} respectively.

The map $\psi$ of Figure~\ref{Xoverline} is obtained by the universal property of the pushout. Thus, it is bijective on states. It then suffices to prove that the map \[\PA^1\psi:\liminj \D^f(1) \to \PA^1\mathbb{M}^\mathcal{P}(X) = \PA^{top}X\] is a homeomorphism since $\mathcal{G}\subset \mathcal{P}$. By Theorem~\ref{calcul_final_structure}, every execution path of $X$ can be written as a finite Moore composition $(f_1\gamma_1\mu_{\ell_1}) * \dots * (f_n\gamma_n\mu_{\ell_n})$ with $n\geq 1$ such that $\ell_1 + \dots + \ell_n = 1$ and such that $f_i = f$ and $\gamma_i$ is an execution path of $A$ or $f_i=\widehat{g}$ and $\gamma_i=\delta_{z_i}\phi_i$ with $z_i\in \mathbf{D}^{n}\backslash\mathbf{S}^{n-1}$ and some $\phi_i\in \mathcal{P}(1,1)$. Let $\overline{f_i}= \overline{g}$ if $f_i=\widehat{g}$ and $\overline{f_i}= \overline{f}$ if $f_i=f$ for $i\in \{1,\dots,n\}$. This gives rise to the execution path $\PA\overline{f_1}(\gamma_1\mu_{\ell_1}) * \dots * \PA\overline{f_n}(\gamma_n\mu_{\ell_n})$ of the Moore flow $\overline{X}$. By the commutativity of the diagram of Figure~\ref{Xoverline}, we obtain the equality 
\[
(f_1\gamma_1\mu_{\ell_1}) * \dots * (f_n\gamma_n\mu_{\ell_n}) = (\PA^1 \psi)\bigg(\PA\overline{f_1}(\gamma_1\mu_{\ell_1}) * \dots * \PA\overline{f_n}(\gamma_n\mu_{\ell_n})\bigg).
\]
This means that the map of Moore flows $\psi:\overline{X} \to \moore^{\mathcal{P}}(X)$ induces a surjective continuous map from $\PA^1\overline{X}$ to $\PA^{top}X$. In other terms, the map $\PA^1\psi$ is a surjection.

Let $\underline{n} = (u_0,\epsilon_1,u_1),(u_1,\epsilon_2,u_2),\dots ,(u_{n-1},\epsilon_n,u_n))$. By \cite[Corollary~5.13]{Moore1}, the topological space $(Z_{u_0,u_1}\ot Z_{u_1,u_2} \ot \dots \ot Z_{u_{n-1},u_n})(1)$ is the quotient of \[\coprod_{(\ell_1,\dots,\ell_n)} \mathcal{P}(1,\ell_1+\dots+\ell_n) \p Z_{u_0,u_1}(\ell_1) \p \dots \p Z_{u_{n-1},u_n}(\ell_n)\] by the equivalence relation generated by the identifications \[(\phi,x_1\phi_1,\dots,x_p\phi_p) \sim ((\phi_1\ot \dots \ot \phi_p)\phi,x_1,\dots,x_p)\] for $\phi\in \mathcal{P}(1,\ell_1+\dots+\ell_p)$, $\phi_i\in \mathcal{P}(\ell_i,\ell'_i)$ and $x_i\in Z_{u_{i-1},u_i}(\ell'_i)$ for $1\leqslant i\leqslant n$. Assume that \[(\phi,\gamma_1,\dots,\gamma_n) \in \mathcal{P}(1,\ell_1+\dots+\ell_n) \p Z_{u_0,u_1}(\ell_1) \p \dots \p Z_{u_{n-1},u_n}(\ell_n)\] is a representative of $\underline{x}$ in $\D^f(\underline{n})$ with $\underline{x}$ simplified. Then \[\PA^1\psi(\underline{x})=\big((f_1\gamma_1)* \dots (f_n\gamma_n)\big)\phi\] with $f_i=f$ if $\epsilon_i=0$ and $f_i=\widehat{g}$ if $\epsilon_i=1$. Using Proposition~\ref{decomposition-tenseur}, write $\phi=\phi_1\ot \dots \ot \phi_n$ with $\phi_i:\ell'_i\to \ell_i$ for $1\leqslant i\leqslant n$ for some $\ell'_1,\dots,\ell'_n$ such that $\ell'_1 + \dots + \ell'_n=1$. Then one has \[(\phi,\gamma_1,\dots,\gamma_n)\sim (\id_1,\gamma_1\phi_1,\dots,\gamma_n\phi_n)\] in $\D^f(\underline{n})$ and therefore \[\PA^1\psi(\underline{x})=(f_1\gamma_1\phi_1)* \dots (f_n\gamma_n\phi_n).\] Recall that \[d((u_0,\epsilon_1,u_1),(u_1,\epsilon_2,u_2),\dots ,(u_{n-1},\epsilon_n,u_n)) = n + \sum_i \epsilon_i.\] As already explained in the proof of \cite[Theorem~7.2]{Moore2} with a lot of details in the case $\mathcal{P}=\mathcal{G}$, since $\underline{x}$ is simplified by hypothesis, it is impossible to have $\epsilon_i=\epsilon_{i+1}=0$ for some $1\leqslant i <n$. There is a composition map starting from $\underline{n}$ in $\mathcal{P}^{g(0),g(1)}(A^0)$ otherwise, which identifies $\underline{x}$ to some $\underline{y}\in \mathcal{D}^f(\underline{m})$ in $\liminj \mathcal{D}^f(1)$ with $d(\underline{m}) < d(\underline{n})$, and it is a contradiction. As also already seen in the proof of \cite[Theorem~7.2]{Moore2} in the case $\mathcal{P}=\mathcal{G}$, if $\epsilon_i=1$, then $\gamma_i = \overline{g}\delta_{z_i}\psi_i$ with $z_i\in \mathbf{D}^n\backslash \mathbf{S}^{n-1}$ and $\psi_i\in \mathcal{P}(\ell_i,1)$. Indeed, if $z_i\in \mathbf{S}^{n-1}$, then there is an inclusion map whose image contains $\underline{x}$, which means also that $\underline{x}$ is identified to some $\underline{y}\in \mathcal{D}^f(\underline{m})$ in $\liminj \mathcal{D}^f(1)$ with $d(\underline{m}) < d(\underline{n})$, which contradicts the fact that $\underline{x}$ is simplified. 

This means that the finite Moore composition $(f_1\gamma_1\phi_1)* \dots (f_n\gamma_n\phi_n)$ is one of the finite Moore compositions given by Theorem~\ref{calcul_final_structure}. Consider another simplified element $\underline{x}'$ in $\D^f(\underline{n}')$ such that $\PA^1\psi(\underline{x}) = \PA^1\psi(\underline{x}')$. It gives rise to another finite Moore composition $(f'_1\gamma'_1\phi'_1)* \dots (f'_{n'}\gamma'_n\phi_{n'})$ as the ones given by Theorem~\ref{calcul_final_structure}. Using Theorem~\ref{calcul_final_structure}, we deduce that $n=n'$ and that
\begin{equation}
	\label{reg1} \tag{R}
	\forall 1\leqslant i \leqslant n, r_i=\natgl(\gamma_i)=\natgl(\gamma'_i), \gamma_i=r_i\eta_i, \gamma'_i=r_i\eta'_i
\end{equation}
and 
\begin{equation}
	\label{reg2} \tag{P}
	(\eta_1\phi_1)\ot \dots \ot (\eta_n\phi_n) = (\eta'_1\phi'_1)\ot \dots \ot (\eta'_n\phi'_n).
\end{equation} We then obtain in $(Z_{u_0,u_1}\ot Z_{u_1,u_2} \ot \dots \ot Z_{u_{n-1},u_n})(1)$ the following sequence of identifications:
\begin{align*}
	(\phi,\gamma_1,\dots,\gamma_n) & \sim ((\eta_1\ot ... \ot \eta_n)\phi,r_1,\dots,r_n)   \\
	& = ((\eta'_1\ot ... \ot \eta'_n)\phi',r_1,\dots,r_n)   \\
	& \sim  (\phi',\gamma'_1,\dots,\gamma'_n),
\end{align*}
the first and third identifications by \eqref{reg1} and \cite[Corollary~5.13]{Moore1} and the equality by \eqref{reg2}. This means that $\underline{x}=\underline{x}'$ in $\liminj \mathcal{D}^f(1)$ and, therefore, that the map $\PA^1\psi$ is one-to-one.

At this point, it is proved that the map $\PA^1\psi:\liminj \D^f(1) \to \PA^{top}X$ is a continuous bijection with $\liminj \D^f(1)$ equipped with the final topology. When we work with the category of $\Delta$-Hausdorff $\Delta$-generated spaces, we deduce that $\liminj \D^f(1)$ equipped with the final topology is $\Delta$-Hausdorff as well, the space $\PA^{top}X$ being $\Delta$-Hausdorff. So whether we work with $\Delta$-Hausdorff or not $\Delta$-generated spaces, the topology of $\liminj \D^f(1)$ is always the final topology.

By \cite[Corollary~2.3]{Moore2}, we must now prove that for all set maps $\xi : [0,1] \to \liminj \D^f(1)$, if the composite map $(\PA^1\psi)\xi:[0,1]\to \PA^{top}X$ is continuous, then the set map $\xi : [0,1] \to \liminj \D^f(1)$ is continuous as well. By Theorem~\ref{carrier-finite-on-compact}, the set of carriers \[\mathcal{T}=\{\carrier((\PA^1\psi)\xi(u)) \mid u\in [0,1]\}\] is finite. For each carrier $\underline{c}\in \mathcal{T}$, let \[
U_{\underline{c}} = \{u\in [0,1]\mid\carrier(\overline{\xi}(u)) = \underline{c}\}.\] Consider the closure $\widehat{U_{\underline{c}}}$ of $U_{\underline{c}}$ in $[0,1]$. We obtain a finite covering of $[0,1]$ by the closed subsets $\widehat{U_{\underline{c}}}$ for $\underline{c}$ running over $\mathcal{T}$. Each $\widehat{U_{\underline{c}}}$ is compact, metrizable and therefore sequential. Note that $\widehat{U_{\underline{c}}}$ has no reason to be $\Delta$-generated: it could be e.g. the Cantor set which is not $\Delta$-generated because it is not homeomorphic to the disjoint sum of its path-connected components. Fix the carrier $\underline{c}$.  

The end of the proof is the \textit{third reduction} and \textit{sequential continuity} sections of the proof of \cite[Theorem~7.3]{Moore2} with the use of Theorem~\ref{img-closed} instead of \cite[Theorem~6.11]{Moore2}. The argument is sketched for the ease of the reader. It suffices to prove that the restriction \[\xi:\widehat{U_{\underline{c}}}\longrightarrow \liminj \D^f(1)\] is sequentially continuous to complete the proof. Let $(u_n)_{n\geq 0}$ be a sequence of $\widehat{U_{\underline{c}}}$ which converges to $u_\infty$. Then the sequence of execution paths $(\PA^1\psi(\xi(u_n)))_{n\geq 0}$ converges to $\PA^1\psi(\xi(u_\infty))$, and therefore, it converges pointwise. All execution paths $\PA^1\psi(\xi(u_n))$ for $n\geq 0$ and $\PA^1\psi(\xi(u_\infty))$ belong to the image of $\PA^{top}\widehat{g_{\underline{c}}}$ (see Notation~\ref{Xc}), this image being closed in $\PA^{top}X_\lambda$ by Theorem~\ref{img-closed}. Besides, each subsequence of $(\PA^1\psi(\xi(u_n)))_{n\geq 0}$ has a limit point by Theorem~\ref{img-closed}. This limit point is unique since \[\PA^1\psi:\liminj \D^f(1) \longrightarrow \PA^{top}X\] is a bijection. The proof is complete thanks to Lemma~\ref{converge}.
\epf

\begin{cor}  \label{calculation-pathspace}
	Suppose that $A$ is a cellular $\mathcal{P}$-multipointed $d$-space. Consider a pushout diagram of $\mathcal{P}$-multipointed $d$-spaces 
	\[
	\begin{tikzcd}[row sep=3em, column sep=3em]
		\globP(\mathbf{S}^{n-1}) \arrow[d] \arrow[r] & A \arrow[d] \\
		\globP(\mathbf{D}^{n}) \arrow[r] & \cocartesian X
    \end{tikzcd}
	\]
	with $n\geq 0$. Then there is the pushout diagram of Moore flows
	\[
	\begin{tikzcd}[row sep=3em, column sep=3em]
		\moore^{\mathcal{P}}(\globP(\mathbf{S}^{n-1}))=\glob(\mathbb{F}^{\mathcal{P}^{op}}_1\mathbf{S}^{n-1}) \arrow[d] \arrow[r] & \moore^{\mathcal{P}}(A) \arrow[d] \\
		\moore^{\mathcal{P}}(\globP(\mathbf{D}^{n}))=\glob(\mathbb{F}^{\mathcal{P}^{op}}_1\mathbf{D}^{n}) \arrow[r] & \cocartesian \moore^{\mathcal{P}}(X).
	\end{tikzcd}
	\]
\end{cor}

\begin{cor} (\cite[Corollary~7.5]{Moore2} for $\mathcal{G}$ and $\mathcal{M}$) \label{looks-like-leftadjoint}
	Let $X$ be a q-cofibrant $\mathcal{P}$-multipointed $d$-space. Then the $\mathcal{P}$-flow $\moore^{\mathcal{P}}(X)$ is q-cofibrant.
\end{cor}

\bpf
As for \cite[Corollary~7.5]{Moore2}, it is a consequence of Theorem~\ref{cof-accessible} and Corollary~\ref{calculation-pathspace}.
\epf

\bth \label{final}
Consider the adjunction $\lmoore^{\mathcal{P}}\dashv\moore^{\mathcal{P}}$ between $\mathcal{P}$-multipointed $d$-spaces and $\mathcal{P}$-flows. Then the unit map and the counit map induce isomorphisms on q-cofibrant objects. This adjunction is a Quillen equivalence between the q-model structures of $\mathcal{P}$-multipointed $d$-spaces and of $\mathcal{P}$-flows.
\eth

\bpf
From Corollary~\ref{calculation-pathspace} and Theorem~\ref{cof-accessible}, we deduce that the unit map and the counit map are isomorphisms on cellular objects, and then, on q-cofibrant objects since the retract of an isomorphism is an isomorphism. From this fact and the fact that all objects are q-fibrant, we deduce that the Quillen adjunction is a Quillen equivalence. See the proofs of \cite[Theorem~7.6, Corollary~7.9 and Theorem~8.1]{Moore2} for further details.
\epf

\bd \label{def:flow}
	The category of small topologically enriched semicategories is isomorphic to the category of $\mathbf{1}$-flows ($\mathbf{1}$ being the terminal category viewed as a reparametrization category). This category is denoted by $\dtop$ and its objects are called \textit{flows} (without using the prefixes Moore or $\mathbf{1}$).
\ed

\begin{nota}
	Let $\C$ be a small category. The \textit{constant diagram functor} is denoted by \[\Delta_{\C}:\C\longrightarrow \topspace.\] 
\end{nota}

By \cite[Proposition~10.5]{Moore1}, the constant diagonal functor induces a functor denoted by \[\moore:\dtop\longrightarrow \dtopP\] such that $\moore(X)^0=X^0$ and such that $\PA_{\alpha,\beta}\moore(X)= \Delta_{\mathcal{P}^{op}}(\PA_{\alpha,\beta}X)$ for all $(\alpha,\beta)\in X^0\p X^0$. By \cite[Proposition~10.6]{Moore1}, for any $\mathcal{P}$-flow $Y$, the data 
\begin{itemize}[leftmargin=*]
	\item The set of states is $Y^0$
	\item For all $\alpha,\beta\in Y^0$, let $Y_{\alpha,\beta}=\liminj \PA_{\alpha,\beta}Y$
	\item For all $\alpha,\beta,\gamma\in Y^0$, the composition law $Y_{\alpha,\beta}\p Y_{\beta,\gamma}\to Y_{\alpha,\gamma}$
\end{itemize}
assemble to a flow denoted by $\lmoore(Y)$. It yields a well-defined functor \[\lmoore:\dtopP \to \dtop.\]
By \cite[Theorem~10.9]{Moore1}, the pair of functors $(\lmoore,\moore)$ gives rise to a Quillen equivalence \[\lmoore \dashv \moore.\]

\begin{nota}
	Let $\dcat = \lmoore \moore^{\mathcal{P}}: \ptop{\mathcal{P}}\longrightarrow \dtop$. 
\end{nota}

\begin{nota} 
	Let $X$ be a $\mathcal{P}$-multipointed $d$-space. Let $\alpha,\beta\in X^0$. The identity of $\lmoore \moore^{\mathcal{P}}(X)$ gives rise to a map $\moore^{\mathcal{P}}(X) \to \moore \lmoore \moore^{\mathcal{P}}(X)$, and by applying $\PA_{\alpha,\beta}^1$ to a continuous map $[-]_{\alpha,\beta} : \PA_{\alpha,\beta}^{top}X \to \PA_{\alpha,\beta} \dcat(X)$.
\end{nota}

\bth  \label{GM-cat} (\cite[Theorem~8.11]{Moore2} and \cite[Theorem~8.14]{Moore2} for $\mathcal{G}$ and $\mathcal{M}$)
Let $X$ be a $\mathcal{P}$-multipointed $d$-space. Then one has (cf. Notation~\ref{trace-space})
\[\dcat(X)^0=X^0\hbox{ and } \forall (\alpha,\beta)\in X^0\p X^0, \PA_{\alpha,\beta}\dcat(X)= \PA_{\alpha,\beta} X\]
and the canonical map $[-]_{\alpha,\beta} : \PA_{\alpha,\beta}^{top}X \to \PA_{\alpha,\beta} \dcat(X)$ is the quotient map which takes an execution path of $X$ from $\alpha$ to $\beta$ to its equivalence class up to reparametrization by $\mathcal{P}(1,1)$. The functor $\dcat: \ptop{\mathcal{P}}\longrightarrow \dtop$ takes q-cofibrant $\mathcal{P}$-multipointed $d$-spaces to q-cofibrant flows. Its total left derived functor in the sense of \cite{HomotopicalCategory} induces an equivalence of categories between the homotopy categories of the q-model structures. 
\eth

\bpf
The equality $\PA_{\alpha,\beta}\dcat(X)= \PA_{\alpha,\beta} X$ comes from the definition of the colimit: see the proof of \cite[Theorem~8.11]{Moore2}. The rest of the proof is mutatis mutandis the proof of \cite[Theorem~8.14]{Moore2} by replacing $\mathcal{G}$ by $\mathcal{P}$ and by using Theorem~\ref{final}. We recall the definition of the functors for the convenience of the reader: 
\[
\begin{tikzcd}[row sep=0.5em, column sep=2.5em]
	(\mathbf{L}\dcat):\ptop{\mathcal{P}}\arrow[r,"(-)^{cof}"]& \ptop{\mathcal{P}} \arrow[r,"\dcat"]& \dtop \\
	(\mathbf{L}\dcat)^{-1}:\dtop \arrow[r,"\moore"] & \dtopP \arrow[r,"(-)^{cof}"]& \dtopP \arrow[r,"\lmoore^{\mathcal{P}}"]& \ptop{\mathcal{P}}
\end{tikzcd}
\]
where $(-)^{cof}$ is a q-cofibrant replacement functor for the corresponding category. 
\epf

\cite[Theorem IV.3.10 and Theorem~IV.3.14]{model2} state in the language of this paper that, for all cellular $\mathcal{G}$-multipointed $d$-spaces $X$ and all $\alpha,\beta\in X^0$, the quotient map $\PA_{\alpha,\beta}^{top}X \to \PA_{\alpha,\beta} X$ is a trivial h-fibration of $\topspace$. The language used in \cite{model2} is different because $\mathcal{G}$-multipointed $d$-spaces were introduced four years later in \cite{mdtop}. 

Consequently, the quotient map $\PA_{\alpha,\beta}^{top}X \to \PA_{\alpha,\beta} X$ is also a trivial h-fibration of $\topspace$ for all q-cofibrant $\mathcal{G}$-multipointed $d$-spaces $X$ and all $\alpha,\beta\in X^0$, the retract of a trivial h-fibration being a trivial h-fibration. The proof relies on \cite[Theorem~III.5.2]{model2} which is fixed in \cite[Theorem~6.8]{leftproperflow}. The latter states that the quotient map $\PA_{\alpha,\beta}^{top}X \to \PA_{\alpha,\beta} X$ has always a section when $X$ is a cellular $\mathcal{G}$-multipointed $d$-space. 

It is not clear whether this section exists for all cellular $\mathcal{M}$-multipointed $d$-spaces. This means that the proof of \cite{model2} does not seem to be generalizable to the case $\mathcal{P}=\mathcal{M}$. The best that can be said is that the quotient map $\PA_{\alpha,\beta}^{top}X \to \PA_{\alpha,\beta} X$ has a section for the so-called \textit{regular} $\mathcal{M}$-multipointed $d$-spaces (the proof will not be given here because it is off-topic). A regular $\mathcal{M}$-multipointed $d$-space is, by definition, a cellular $\mathcal{M}$-multipointed $d$-space $X_\lambda$ such that for all globular cells $c_\nu$ of $X_\lambda$ and all $z\in \mathbf{D}^{n_\nu}$, the execution path $\widehat{g_{\nu}}\delta_{z}$ is regular. Remember that, in general, the execution paths $\widehat{g_{\nu}}\delta_{z}$ are regular only when $z\in \mathbf{D}^{n_\nu}\backslash \mathbf{S}^{n_\nu-1}$ by Proposition~\ref{example_reg}. In plain English, a regular $\mathcal{M}$-multipointed $d$-space is a cellular $\mathcal{M}$-multipointed $d$-space such that the attaching maps take regular execution paths of the boundary of globular cells to regular execution paths. 

Thanks to the results of this paper, we can still prove that this quotient map is a homotopy equivalence for all cellular $\mathcal{M}$-multipointed $d$-spaces as follows.

\bth \label{quotient-reparam}
Let $X$ be a q-cofibrant $\mathcal{P}$-multipointed $d$-space. Let $\alpha,\beta\in X^0$. Then the quotient map \[\PA_{\alpha,\beta}^{top}X \longrightarrow \PA_{\alpha,\beta} X\] is a homotopy equivalence from an m-cofibrant space to a q-cofibrant space.
\eth

\bpf
By Corollary~\ref{looks-like-leftadjoint}, the $\mathcal{P}$-flow $\moore^{\mathcal{P}}(X)$ is q-cofibrant. By \cite[Theorem~9.11]{Moore1}, we deduce that the $\mathcal{P}$-space $\PA_{\alpha,\beta} \moore^{\mathcal{P}}(X)$ is projective q-cofibrant. Using \cite[Corollary~7.2]{dgrtop}, we deduce that the $\mathcal{P}$-space $\PA_{\alpha,\beta} \moore^{\mathcal{P}}(X)$ is injective m-cofibrant. This implies that the topological space $\PA_{\alpha,\beta}^{top}X = \PA_{\alpha,\beta}^1 \moore^{\mathcal{P}}(X)$ is m-cofibrant. By Theorem~\ref{GM-cat}, the flow $\dcat(X)$ is q-cofibrant. Using \cite[Theorem~9.11]{Moore1} again or \cite[Theorem~5.7]{leftproperflow}, we deduce that the topological space $\PA_{\alpha,\beta}\dcat(X)=\PA_{\alpha,\beta}X$ is q-cofibrant. The $\mathcal{P}$-flow $\moore^{\mathcal{P}}(X)$ being q-cofibrant, from the Quillen equivalence $\lmoore\dashv \moore$, we obtain the weak equivalence of $\mathcal{P}$-flows 
\[
\moore^{\mathcal{P}}(X) \longrightarrow \moore \lmoore \moore^{\mathcal{P}}(X).
\]
By applying $\PA_{\alpha,\beta}^1(-)$ to both sides, we deduce a weak homotopy equivalence of topological spaces 
\[
\PA_{\alpha,\beta}^{top}X = \PA_{\alpha,\beta}^1\moore^{\mathcal{P}}(X)\longrightarrow \PA_{\alpha,\beta}^1\moore \lmoore \moore^{\mathcal{P}}(X) = \PA_{\alpha,\beta}^1\moore (\dcat(X)) = \PA_{\alpha,\beta}X
\]
from an m-cofibrant space to a q-cofibrant space, i.e. between two m-cofibrant spaces. By \cite[Corollary~3.4]{mixed-cole}, we deduce that this map is a homotopy equivalence.
\epf

We do not know whether this quotient map is still, at least, a weak homotopy equivalence for a general $\mathcal{P}$-multipointed $d$-space. It is likely that the saturation hypothesis introduced in Section~\ref{section-saturated} plays a role.

Proposition~\ref{final2} should have been put in \cite{Moore1} as an application of the results of the latter paper: it is an omission. It is used in Theorem~\ref{comparisonGM}. The inclusion functor $i:\mathcal{G}\subset \mathcal{M}$ induces an enriched functor \[i^*:\topdgrM \longrightarrow \topdgrG\] from $\mathcal{M}$-spaces to $\mathcal{G}$-spaces. It is a right adjoint between the underlying categories, the left adjoint being the enriched left Kan extension along $i$ given by the formula 
\[
\lan_i(D) = \int^{\ell} \mathcal{M}(-,i(\ell)) \p D(\ell).
\]

\bp \label{final2}
The functor $i^*:\topdgrM \to \topdgrG$ induces a functor \[i^*:\dtopM \longrightarrow \dtopG\] which is a right Quillen equivalence between the q-model structures of $\dtopM$ and $\dtopG$. 
\ep

\bpf
By \cite[Section~6]{Moore1}, a $\mathcal{P}$-flow consists of a set of states $X^0$, for each pair $(\alpha,\beta)$ of states a $\mathcal{P}$-space $\PA_{\alpha,\beta}X$ of $\topdgrP_0$ and for each triple $(\alpha,\beta,\gamma)$ of states an associative composition law $*:\PA^{\ell_1}_{\alpha,\beta}X \p \PA^{\ell_2}_{\beta,\gamma}X \to \PA^{\ell_1+\ell_2}_{\alpha,\gamma}X$ which is natural with respect to $(\ell_1,\ell_2)$ in an obvious way. From an $\mathcal{M}$-flow $D$, we therefore obtain a $\mathcal{G}$-flow $i^*(D)$ with $D^0=i^*(D)^0$ and $\PA_{\alpha,\beta}i^*(D) = i^*(\PA_{\alpha,\beta}D)$. By the explicit calculation of limits in $\dtopM$ and in $\dtopG$ made in \cite[Theorem~6.8]{Moore1}, and since limits are calculated objectwise in $\topdgrM_0$ and $\topdgrG_0$ by \cite[Proposition~5.3]{dgrtop}, the functor $i^*:\dtopM \to \dtopG$ is limit-preserving. By \cite[Theorem~6.13]{Moore1}, the $\mathcal{P}$-space of execution paths functor $\PA:\dtopP\to \topdgrP_0$ of Definition~\ref{Pmoore} is a right adjoint for any reparametrization category $\mathcal{P}$~\footnote{Note that this fact holds because we work with locally presentable categories: see \cite[Theorem~5.10]{leftproperflow}.}. Therefore it is accessible by \cite[Theorem~1.66]{TheBook}. Since colimits are calculated objectwise in $\topdgrM_0$ and $\topdgrG_0$ by \cite[Proposition~5.3]{dgrtop}, the functor $i^*:\dtopM \to \dtopG$ is then accessible. Therefore it is a right adjoint by \cite[Theorem~1.66]{TheBook}. The functor $i^*:\dtopM \to \dtopG$ preserves q-fibrations and trivial q-fibrations by definition of the q-model structures. Consequently, it is a right Quillen adjoint. Thus the commutative diagram of right adjoints 
\[
\begin{tikzcd}[row sep=3em, column sep=3em]
	\topspace \arrow[r,equal] \arrow[d,"\Delta_{\mathcal{M}^{op}}"'] & \topspace \arrow[d,"\Delta_{\mathcal{G}^{op}}"]\\
	\topdgrM \arrow[r,"i^*"] & \topdgrG
\end{tikzcd}
\]
gives rise by \cite[Proposition~10.7]{Moore1} to the commutative diagram of right Quillen adjoints
\[
\begin{tikzcd}[row sep=3em, column sep=3em]
	\dtop \arrow[r,equal] \arrow[d,"\moore"'] & \dtop \arrow[d,"\moore"]\\
	\dtopM \arrow[r,"i^*"] & \dtopG
\end{tikzcd}
\]
where $\dtop$ is equipped with its q-model structure. By \cite[Theorem~10.9]{Moore1}, the two vertical right Quillen adjoints are right Quillen equivalences. The proof is complete thanks to the \ttt.
\epf

We conclude with the following comparison theorem:

\bth \label{comparisonGM}
The inclusion functor $i:\mathcal{G}\subset \mathcal{M}$ induces a functor \[j:\ptop{\mathcal{M}} \longrightarrow \ptop{\mathcal{G}}.\] There is the commutative square of right Quillen equivalences between the four q-model structures
\[
\begin{tikzcd}[row sep=3em, column sep=3em]
	\ptop{\mathcal{M}} \arrow[r,"j"] \arrow[d,"\moore^{\mathcal{M}}"'] & \ptop{\mathcal{G}} \arrow[d,"\moore^{\mathcal{G}}"]\\
	\dtopM \arrow[r,"i^*"] & \dtopG
\end{tikzcd}
\]
\eth

\bpf
It is easy to see that the diagram is commutative: each functor is a forgetful functor indeed. The forgetful functor $\Omega:\ptop{\mathcal{P}}\to \mtop$ from $\mathcal{P}$-multipointed $d$-spaces to multipointed spaces being topological by Theorem~\ref{final-structure-revisited} for $\mathcal{P}$ equal to $\mathcal{G}$ or $\mathcal{M}$, the functor $j:\ptop{\mathcal{M}} \to \ptop{\mathcal{G}}$ is limit-preserving and finitely accessible: \textit{finitely} because a multipointed $d$-space is equipped with a \textit{set} of execution paths and because the $\Omega$-final structure is given by the \textit{finite} Moore compositions by Theorem~\ref{final-structure-revisited}. By \cite[Theorem~1.66]{TheBook}, the functor $j:\ptop{\mathcal{M}} \to \ptop{\mathcal{G}}$ is therefore a right adjoint. It takes (trivial resp.) q-fibrations to (trivial resp.) q-fibrations by definition of them. Thus it is a right Quillen adjoint. The two vertical functors are right Quillen equivalences by Theorem~\ref{final}. The bottom horizontal functor is a right Quillen equivalence by Proposition~\ref{final2}. The proof is complete thanks to the \ttt.
\epf

\begin{nota}
	Write $\mathcal{F}_{\mathcal{G}}^{\mathcal{M}} : \ptop{\mathcal{G}} \to \ptop{\mathcal{M}}$ for the left adjoint of the inclusion functor $j:\ptop{\mathcal{M}}\subset \ptop{\mathcal{G}}$. 
\end{nota}

The unit of the adjunction \[X\longrightarrow j(\mathcal{F}_{\mathcal{G}}^{\mathcal{M}}(X))\] preserves the underlying space and the set of states. It induces a map from the space of execution paths of $X$ to its closure under the reparametrization by all maps of $\mathcal{M}$. It is a weak homotopy equivalence when $X$ is a q-cofibrant $\mathcal{G}$-multipointed $d$-space by Theorem~\ref{comparisonGM}: this assertion is also a consequence of Corollary~\ref{calculation-pathspace} and Theorem~\ref{cof-accessible} and of the fact that \[\mathcal{F}_{\mathcal{G}}^{\mathcal{M}}(\globG(Z))=\globM(Z)\] for all topological spaces $Z$. The counit map \[\mathcal{F}_{\mathcal{G}}^{\mathcal{M}}(j(Y)) \stackrel{\iso}\longrightarrow Y\] is an isomorphism for all $\mathcal{M}$-multipointed $d$-spaces $Y$ by definition of $\mathcal{F}_{\mathcal{G}}^{\mathcal{M}}$. By Theorem~\ref{comparisonGM}, we deduce that $\mathcal{F}_{\mathcal{G}}^{\mathcal{M}}(j(Y)^{cof})$ is a q-cofibrant replacement of $Y$ in $\ptop{\mathcal{M}}$ where $j(Y)^{cof}$ is a q-cofibrant replacement of $j(Y)$ in $\ptop{\mathcal{G}}$. The latter fact can be proved directly by obtaining a q-cofibrant replacement by the small object argument.

\section{The saturation hypothesis}
\label{section-saturated}

This notion is very important in DAT. It appears in various forms in the literature: see \cite[Definition~4.3]{reparam} \cite[Remark~4.3]{FR}  \cite[Definition~2.9]{zbMATH06404305} \cite[Definition~2.18]{Ziemiaski2012}. In \cite{Ziemiaski2012}, this notion leads to an isomorphism between two categories of continuous (i.e. non-multipointed) geometric models of concurrency, namely the full subcategory of saturated $d$-spaces of the category of Grandis' $d$-spaces in the sense of \cite{mg} and the full subcategory of good streams of the category of Krishnan's streams in the sense of \cite{MR2545830}. It is worth nothing that all examples coming from concurrency theory are saturated.

The saturation hypothesis is a closure property. In the setting of this paper, the idea is that a continuous path of the underlying space of a $\mathcal{P}$-multipointed $d$-space such that there exists a reparametrization which is an execution path should be an execution path as well. This way, we avoid pathological behaviors like the one of the multipointed $d$-space $\mathbb{I}_\phi$ defined in the proof of Proposition~\ref{ex-non-saturated}. All $\mathcal{G}$-multipointed $d$-spaces are automatically saturated because all maps of $\mathcal{G}$ are invertible. The purpose of this section is to prove that the saturation hypothesis can be added to the definition of an $\mathcal{M}$-multipointed $d$-space without changing the main results of this paper.

\begin{nota}
	Let $X$ be a $\mathcal{P}$-multipointed $d$-space. Consider for all $\alpha,\beta\in X^0$ the set of continuous paths 
	\[
	\widehat{\PA}^{top}_{\alpha,\beta} X = \{\gamma\in \topspace([0,1],|X|)\mid \exists \phi \in \mathcal{P}(1,1),\gamma\phi\in \PA^{top}_{\alpha,\beta}X\}.
	\]
	Let \[\widehat{\PA}^{top}X = \coprod_{(\alpha,\beta)\in X^0\p X^0} \widehat{\PA}^{top}_{\alpha,\beta} X.\] 
\end{nota}

\bp
Let $X$ be a $\mathcal{P}$-multipointed $d$-space. The triple $\widehat{X}=(|X|,X^0,{\PA}^{top}\widehat{X})$ with ${\PA}^{top}\widehat{X}=\widehat{\PA}^{top}X$ is a $\mathcal{P}$-multipointed $d$-space. The inclusion $X\subset \widehat{X}$ yields a natural map of multipointed $d$-spaces.
\ep

\bpf
When $\mathcal{P}=\mathcal{G}$, one has $\widehat{X}=X$ since all maps of $\mathcal{G}$ are invertible. 

Assume for the rest of the proof that $\mathcal{P}=\mathcal{M}$. Let $(\gamma_1,\gamma_2)\in \widehat{\PA}^{top}_{\alpha,\beta} X \p \widehat{\PA}^{top}_{\beta,\gamma} X$. Then there exist $\phi_1,\phi_2\in \mathcal{M}(1,1)$ such that $\gamma_1\phi_1\in {\PA}^{top}_{\alpha,\beta} X$ and $\gamma_2\phi_2\in {\PA}^{top}_{\beta,\gamma} X$. Write 
\begin{align*}
	(\gamma_1\phi_1) *_N (\gamma_2\phi_2) & = (\gamma_1\phi_1\mu_{\frac{1}{2}}) * (\gamma_2\phi_2\mu_{\frac{1}{2}}) \\
	& = (\gamma_1 * \gamma_2) \bigg((\phi_1\mu_{\frac{1}{2}}) \ot (\phi_2\mu_{\frac{1}{2}})\bigg) \\
	& = \bigg((\gamma_1\mu_{\frac{1}{2}}) * (\gamma_2\mu_{\frac{1}{2}})\bigg)\bigg(\mu_{\frac{1}{2}}^{-1}\ot \mu_{\frac{1}{2}}^{-1}\bigg) \bigg((\phi_1\mu_{\frac{1}{2}}) \ot (\phi_2\mu_{\frac{1}{2}})\bigg)\\
	& = (\gamma_1 *_N \gamma_2) \bigg((\mu_{\frac{1}{2}}^{-1}\phi_1\mu_{\frac{1}{2}}) \ot (\mu_{\frac{1}{2}}^{-1}\phi_2\mu_{\frac{1}{2}})\bigg),
\end{align*}
the first and last equalities by definition of the normalized composition, the two other equalities by definition of the Moore composition of paths and by definition of $\otimes$. We deduce that $\gamma_1 *_N \gamma_2\in \widehat{\PA}^{top}_{\alpha,\gamma} X$. 

Let $\gamma\in \widehat{\PA}^{top}_{\alpha,\beta} X$. Then there exists $\phi\in \mathcal{M}(1,1)$ such that $\gamma\phi\in {\PA}^{top}_{\alpha,\beta} X$ . Let $\psi\in \mathcal{M}(1,1)$. Then by \cite[Proposition~2.19]{reparam}, there exist $\xi_1,\xi_2\in \mathcal{M}(1,1)$ such that $\psi\phi\xi_1 = \phi \xi_2$: geometrically, the two execution paths  $\psi\phi$ and $\phi$ of $\vI^{\mathcal{M}}$ are reparametrization equivalent. We obtain $\gamma\psi\phi\xi_1 = \gamma \phi \xi_2$. From $\gamma\phi\in {\PA}^{top}_{\alpha,\beta} X$, we deduce that $\gamma\psi\phi\xi_1 \in {\PA}^{top}_{\alpha,\beta} X$ as well, and therefore that $\gamma\psi\in \widehat{\PA}^{top}_{\alpha,\beta} X$. The proof is complete.
\epf

\bd
A $\mathcal{P}$-multipointed $d$-space $X$ is \textit{saturated} if the natural map $X\to \widehat{X}$ is an isomorphism. The full subcategory of saturated $d$-spaces is denoted by $\ptopsat{\mathcal{P}}$. 
\ed

\bp \label{ex-non-saturated}
All $\mathcal{G}$-multipointed $d$-spaces are saturated. There exists an $\mathcal{M}$-multi\-pointed $d$-space which is not saturated. An $\mathcal{M}$-multipoin\-ted $d$-space is saturated if and only if any continuous path which is reparametrization equivalent (see Definition~\ref{reparam-equivalence} and \cite[Definition~1.2]{reparam}) to an execution path is an execution path.
\ep

\bpf
The first and last assertions are clear. Let $\phi\in \mathcal{M}(1,1)\backslash \mathcal{G}(1,1)$. Consider the $\mathcal{M}$-multipointed $d$-space $\mathbb{I}_\phi$ defined by $(|\mathbb{I}_\phi|,\mathbb{I}_\phi^0) = ([0,1],\{0,1\})$ and $\PA^{top}\mathbb{I}_\phi= \PA^{top}_{0,1}\mathbb{I}_\phi = \{\phi\psi \mid \psi\in \mathcal{P}(1,1)\}$. Then the $\mathcal{M}$-multipointed $d$-space $\mathbb{I}_\phi$ is not saturated. 
\epf

\bp \label{saturated-locally-presentable}
The category $\ptopsat{\mathcal{M}}$ is a reflective locally presentable full subcategory of $\ptop{\mathcal{M}}$.
\ep

\bpf
The left adjoint of the inclusion $\ptopsat{\mathcal{M}} \subset \ptop{\mathcal{M}}$ is given by the functor $X\mapsto \widehat{X}$. The category $\ptopsat{\mathcal{M}}$ is axiomatized by the theory described in the proof of Proposition~\ref{locally-presentable} with the additional axioms $(\forall x) R(x.t) \Rightarrow R(x)$ with $t\in \mathcal{M}(1,1)$ (see \cite[Remark~4.3]{FR}) where $R$ is the $[0,1]$-ary relational symbol encoding execution paths. By \cite[Theorem~5.30]{TheBook}, the proof is complete.
\epf

\bp \label{saturated-retract}
A retract of a saturated $\mathcal{M}$-multipointed $d$-space is saturated.
\ep

\bpf
Let $X$ be a retract of a saturated $\mathcal{M}$-multipointed $d$-space $Y$. This means that the identity of $X$ factors as a composite $r.i:X\to Y \to X$ ($r.i$ means the composition of $r$ and $i$). Let $\gamma$ be a continuous path of $|X|$ and $\phi\in \mathcal{M}(1,1)$ such that $\gamma.\phi$ is an execution path of $X$. Then $i.\gamma.\phi$ is an execution path of $Y$, the continuous map $i:X\to Y$ being a map of $\mathcal{M}$-multipointed $d$-spaces. Since $Y$ is saturated by hypothesis, we deduce that $i.\gamma$ is an execution path of $Y$. Thus $\gamma=r.i.\gamma$ is an execution path of $X$, the continuous map $r:Y\to X$ being a map of $\mathcal{M}$-multipointed $d$-spaces.
\epf

\bp \label{saturated-cofibrant}
All q-cofibrant $\mathcal{M}$-multipointed $d$-spaces are saturated.
\ep

\bpf
First of all, observe that the $\mathcal{M}$-multipointed $d$-space $\globM(Z)$ is saturated for all topological spaces $Z$. The functor $X\mapsto \widehat{X}$ from $\ptop{\mathcal{M}}$ to $\ptopsat{\mathcal{M}}$ being a left adjoint, we deduce that all cellular $\mathcal{M}$-multipointed $d$-spaces are saturated. The proof is complete thanks to Proposition~\ref{saturated-retract}.
\epf

\bth
The q-model structure of $\ptop{\mathcal{M}}$ can be lifted (in the sense of \cite{HKRS17,GKR18}) along the right adjoint $\ptopsat{\mathcal{M}}\subset \ptop{\mathcal{M}}$. It is a combinatorial model structure, called the q-model structure of $\ptopsat{\mathcal{M}}$, such that:
\begin{itemize}
	\item A set of generating cofibrations is 
	\[\{\globM(\mathbf{S}^{n-1})\subset \globM(\mathbf{D}^{n}) \mid n\geq 0\} \cup \{C:\varnothing \to \{0\},R:\{0,1\} \to \{0\}\},\]
	the maps between globes being induced by the closed inclusions $\mathbf{S}^{n-1}\subset \mathbf{D}^{n}$.
	\item A set of generating trivial cofibrations is \[
	\{\globM(\mathbf{D}^{n})\subset \globM(\mathbf{D}^{n+1}) \mid n\geq 0\},
	\]
	the maps between globes being induced by the closed inclusions $(x_1,\dots,x_n)\mapsto (x_1,\dots,x_n,0)$.
	\item The weak equivalences are the maps of saturated $\mathcal{M}$-multipointed $d$-spaces $f:X\to Y$  inducing a bijection $f^0:X^0\iso Y^0$ and a weak homotopy equivalence $\PA^{{top}} f:\PA^{{top}} X \to \PA^{{top}} Y$.
	\item The fibrations are the maps of saturated $\mathcal{M}$-multipointed $d$-spaces $f:X\to Y$  inducing a q-fibration $\PA^{{top}} f:\PA^{{top}} X \to \PA^{{top}} Y$ of topological spaces.
\end{itemize}
This model structure is combinatorial and all its objects are fibrant.
\eth

\bpf
The path object for $\mathcal{M}$-multipointed $d$-spaces is described in \cite[Corollary~6.7]{QHMmodel}. Let $X$ be an $\mathcal{M}$-multipointed $d$-space. The $\mathcal{M}$-multipointed $d$-space $\cocyl(X)$ has the underlying topological space $\ttop([0,1],|X|)$, the set of states $X^0$ identified with the corresponding set of constant maps of $\ttop([0,1],|X|)$, and \[\PA^{top}_{\alpha,\beta}\cocyl(X) = \ttop([0,1],\PA^{top}_{\alpha,\beta}X)\] for all $\alpha,\beta\in X^0$. We deduce immediately that if $X$ is saturated, then $\cocyl(X)$ is saturated as well. Thus, the q-model structure of $\ptop{\mathcal{M}}$ can be lifted along the right adjoint $\ptopsat{\mathcal{M}}\subset \ptop{\mathcal{M}}$ by using the Quillen Path Object argument: see \cite{HKRS17,GKR18} and for a presentation in this particular context \cite[Theorem~2.1]{QHMmodel}. The characterizations of the generating cofibrations and generating trivial cofibrations is a consequence of Proposition~\ref{saturated-cofibrant} and of a theorem due to Kan (e.g. see \cite[Theorem~11.3.2]{ref_model2}).
\epf

\bth \label{Quillen-equiv-sat}
The right Quillen adjoint \[\ptopsat{\mathcal{M}}\subset \ptop{\mathcal{M}}\] between the q-model structures is a right Quillen equivalence such that the unit and the counit of the adjunction are isomorphisms on q-cofibrant objects.
\eth

\bpf The q-cofibrant objects of the q-model structures of $\ptopsat{\mathcal{M}}$ and $\ptop{\mathcal{M}}$ are the same. All objects are fibrant. Hence the proof is complete.
\epf

\begin{cor} \label{Quillen-equiv-sat-2}
	The composite functor \[\ptopsat{\mathcal{M}}\subset \ptop{\mathcal{M}} \stackrel{\moore^{\mathcal{M}}}\longrightarrow \dtopM\] from $\mathcal{M}$-multipointed $d$-spaces to Moore flows is a right Quillen equivalence such that the unit and the counit of the adjunction are isomorphisms on q-cofibrant objects.
\end{cor}

\bpf
It is a consequence of Theorem~\ref{final} and Theorem~\ref{Quillen-equiv-sat}.
\epf

\begin{cor}  \label{Quillen-equiv-sat-3}  
	The mapping $X \mapsto \dcat(X)$ induces a functor from $\ptopsat{\mathcal{P}}$ to $\dtop$. It takes q-cofibrant saturated $\mathcal{P}$-multipointed $d$-spaces to q-cofibrant flows. Its total left derived functor in the sense of \cite{HomotopicalCategory} induces an equivalence of categories between the homotopy categories of the q-model structures. 
\end{cor}

\bpf
It is a consequence of Theorem~\ref{GM-cat} and Theorem~\ref{Quillen-equiv-sat}.
\epf

The reason for introducing the saturation hypothesis is to rule out pathological behaviors which are meaningless from a computer scientific viewpoint. There is another reason based on the following conjecture which was already stated for $\mathcal{G}$ in \cite[Conjecture~6.6]{mdtop}.

\begin{conj}
	Take for $\topspace$ the category of $\Delta$-Hausdorff $\Delta$-generated spaces. The q-model structures of $\ptop{\mathcal{G}}$ and $\ptopsat{\mathcal{M}}$ are left proper. 
\end{conj}

Indeed, we suspect that the saturation hypothesis and the $\Delta$-Hausdorff condition are necessary for left properness.


\end{document}